\documentclass[reqno]{amsart}

\usepackage{amsmath}
\usepackage{amsthm}
\usepackage{amssymb}
\usepackage{url}

\usepackage[english]{babel}
\usepackage{amsmath,amsthm,amssymb,amscd,amsfonts,mathrsfs,tikz,bm}
\usepackage{xcolor}
\usepackage[latin1]{inputenc}
\usepackage{soul}
\usepackage[all]{xy}
\usetikzlibrary{matrix,arrows,positioning,calc}
\usepackage[T1]{fontenc}
\usepackage{textcomp}
\usepackage{palatino,helvet}
\usepackage{titlesec}
\usepackage[hang,flushmargin]{footmisc}
\usepackage{enumerate} 
\usepackage{graphicx}

\usepackage{tikz}
\usepackage{tikz-cd}

\titlelabel{\thetitle.\quad}

\titleformat{\subsection}[runin]{\normalfont\bfseries}{\thesubsection.}{3pt}{}
\titleformat{\subsubsection}[runin]{\normalfont\bfseries}{\thesubsubsection.}{3pt}{}

\numberwithin{equation}{subsection}\theoremstyle{plain}

\theoremstyle{definition}

\theoremstyle{remark}

\def\Z{\mathbb Z}

\def\mcE{\mathcal{E}}
\def\mcC{\mathcal{C}}

\def\Qco{\mathfrak Qco}

\def\Mod{\mbox{-Mod}}
\def\St{\operatorname{St}}

\def\bfC+{\mathbf{C}_+}

\def\C{{\mathcal{C}}}
\def\V{{\mathcal{V}}}

\def\A{\mathcal{A}}
\def\bfC{\mathbf{C}}

\def\Qco{\mathfrak{Qcoh}}
\def\V{\mathcal{V}}
\def\P{\mathcal{P}}
\def\A{\mathcal{A}}

\newcommand{\tHom}{\mathrm{Hom}_k^{\scriptscriptstyle\bullet}}
\newcommand{\tTen}[1][k]{\otimes_{#1}^{\scriptscriptstyle\bullet}}
\newcommand{\mHom}{\underline{\mathrm{Hom}}_{\mspace{2mu}k}^{\mspace{1mu}\scriptscriptstyle\bullet}}
\newcommand{\mTen}[1][k]{\underline{\otimes}_{\mspace{2mu}#1}^{\mspace{1mu}\scriptscriptstyle\bullet}}
\newcommand{\shHomqce}{\mathscr{H}\mspace{-2.5mu}om^\mathrm{qc}}
\newcommand{\shHome}{\mathscr{H}\mspace{-2.5mu}om}

\newcommand{\End}{\mbox{End}}

\newcommand{\Add}{\mbox{Add}}
\newcommand{\im}{\operatorname{Im}}

\newcommand{\Hom}{\operatorname{Hom}}

\newcommand{\id}{\operatorname{id}}

\newcommand{\Ker}{\mbox{Ker}}
\newcommand{\coker}{\mbox{Coker}}

\newcommand{\Ab}{\operatorname{Ab}}

\newcommand{\ev}{\operatorname{ev}}

\newcommand{\Qcoh}[1]{\mathsf{Qcoh}(#1)}

\begin{document}

\title[]{Homotopy in Enriched Functor Categories: Frobenius Exact Structures and Weak Serre Functors
}

\author[S. Odaba\c si]{Sinem Odaba\c si}
\address{
Universidad de Murcia \\
Facultad de Matem\'aticas \\
Spain
}
\email[S. Odaba\c si]{sinem.odabasi@um.es}

\date{}

\maketitle

\begin{abstract}
Let $(\V, \otimes, I,[-,-])$ be a closed symmetric monoidal category with finite (co)limits. Assume that   $\V$ is weakly idempotent complete and additive. For a given small $\V$-category $\A$, we investigate sufficient conditions on $\A$ under which the ordinary category $[\A, \V]_0$ of $\V$-functors from $\A$ to $\V$, equipped with the exact structure $\mathcal{E}_{\textrm{lowt}}$, is a Frobenius exact category,and  hence, induces a triangulated stable category. The exact structure $\mathcal{E}_{\textrm{lowt}}$ encompasses  the degre-wise trivial exact structure on  categories of $N$-chain complexes and $m$-periodic chain complexes, with terms in a variety of settings, including left $R$-modules, quasi-coherent sheaves on a scheme, and Banach spaces. We subsequently address the question of describing abstract left and right homotopy relations in $([\A, \V]_0;\mathcal{E}_{\textrm{lowt}})$  in terms of families of morphisms defined objectwise, as in the cases of chain-homotopy and relative homotopy relations. To this end, we introduce the notion of left/right-weak Serre $\V$-functors on $\A$, and investigate their relationships with the classical Serre functor relative to a commutative ring $k$.

\end{abstract}

\section*{Introduction}


One of the classical examples of an algebraic triangulated category in Homological Algebra is the homotopy category $\mathbf{K}(R)$ of a ring $R$, which is obtained from the category $\mathbf{C}(R)$ of chain complexes of left $R$-modules by identifying  chain-homotopic maps. Since the equivalence relation of being  chain-homotopic is compatible with the composition and  the abelian group structure  of chain maps, the class of chain maps that are chain-homotopic to zero forms an ideal $\mathcal{I}$ in $\mathbf{C}(R)$. The homotopy category $\mathbf{K}(R)$ is defined as the quotient of the category $\mathbf{C}(R)$  by the ideal $\mathcal{I}$. For an integer $N \geq 2$, analogous chain-homotopy relation and results have been developed for the category $\mathbf{C}_N(R)$ of $N$-complexes of left $R$-modules and the corresponding homotopy category $\mathbf{K}_N(R)$; see \cite{Gil15}.

On the other hand, for an even integer $m \geq 2$ and an abelian category $\mathcal{C}$, the category $\mathbf{C}_{\equiv_m}(\mathcal{C})$ of $m$-cycle complexes over $\mathcal{C}$ is introduced in \cite[\S~7]{PX97}, where the equivalence relation on $m$-cycle complexes, called relative homotopy, is also defined. In \cite[Proposition~7.1]{PX97}, the authors show that   the so-called relative homotopy category $\mathbf{K}_{\equiv_m}(\mathcal{C})$ 
of $m$-cycle complexes over $\mathcal{C}$, obtained as the quotient  category of $\mathbf{C}_{\equiv_m}(\mathcal{C})$ by the ideal of relative homotopic to zero maps, is a triangulated category.

As a general framework, a classical and abstract approach to constructing nontrivial algebraic triangulated categories  consists in identifying a nontrivial Frobenius exact structure $\mathcal{E}$ (see \eqref{Frob}) on a given additive category $\mathcal{C}$; see \cite[\S~3.6]{Keller}. In this setting, the stable category $\St_{\mathcal{P}}(\mathcal{C})$, obtained by taking quotient of $\mathcal{C}$ by the ideal $\mathcal{P}$ of morphisms factoring through a $\mathcal{E}$-projective (or equivalently, $\mathcal{E}$-injective) objects,   is  an additive triangulated category.  This construction also induces an abstract notion of homotopy between morphisms $f,g: A \rightarrow B$ in $\mathcal{C}$: $f$ and $g$ are said to be homotopic if $f-g$ belongs to $\mathcal{P}(A,B)$.

Now, it is well known that both homotopy categories $\mathbf{K}_N(R)$ and $\mathbf{K}_{\equiv_m}(\mathcal{C})$  can be realized as the stable categories associated with a suitable Frobenius exact structure on $\mathbf{C}_N(R)$ and $\mathbf{C}_{\equiv_m}(R)$, respectively. Namely, the exact structure $\mathcal{E}_{\textrm{dwt}}$ on $\mathbf{C}_N(R)$  consisting of those short exact sequences in $\mathbf{C}(R)$ whose evaluation at each $i\in \Z$ is a trivial short exact sequence in $R\Mod$ is a Frobenius exact structure, and  the ideals $\mathcal{I}$ and $\mathcal{P}$ of  chain-homotopic to zero maps and  homotopic to zero maps in the sense of Frobenius exact categories coincide; see \cite[Corollary~3.5, Proposition~4.1]{Gil15}. Consequently, we obtain an equivalence of triangulated categories $\mathbf{K}_N(R) \cong \St_{\mathcal{P}}(\mathbf{C}(R))$. For the analogous construction and result on $\mathbf{C}_{\equiv_m}(R)$, see  \cite[proof of Proposition~7.1-(1)]{PX97}.

Both categories $\mathbf{C}_N(R)$ and $\mathbf{C}_{\equiv_m}(R)$ are known to be particular cases of $R$-module valued additive functor categories   from a small preadditive category $\A$. This naturally raises the following question: Can one define an appropriate  `homotopy category' for $\Add(\mathcal{A},R\Mod)$ for a more general choice of $\mathcal{A}$? Rather than relying on an intrinsic notion of homotopy, as in the preceding examples, we approach this question from the exact category perspective. More precisely, we firstly investigate the sufficient conditions on the domain category $\mathcal{A}$ under which $\Add(\mathcal{A},R\Mod)$ admits a nontrivial Frobenius exact structure covering the degree-wise split exact structures on $\mathbf{C}_N(R)$ and $\mathbf{C}_{\equiv_m}(R)$. When we have it,  we also investigate when the corresponding  abstract notion of homotopy admits a more explicit, intrinsic description, analogous to the chain-homotopy relation in $\mathbf{C}_N(R)$ and the relative homotopy relation in $\mathbf{C}_{\equiv_m}(R)$. These questions will be addressed in Section~4.  The central notion introduced in this paper to achive it is that of a \textit{right/left-weak  Serre functors} on $\mathcal{A}$ (Definition~\ref{def:left_right_weak_Serre}),  which we discuss in greater detail below.

\sloppy In order to encompass homotopy categories arising from certain non-module additive categories  such as the categories of quasi-coherent sheaves over a scheme $X$, $\Qco(X)$, and of Banach spaces, $\mathbf{Ban}$, we propose to investigate those questions within the framework of enriched category theory. The category $\mathbf{Ban}$ is of particular interest because it is not abelian but rather quasi-abelian with  finite (co)limits (Example~\ref{ex:Banach}). This case therefore forces us to carefully verify  and  reformulate certain known results from the literature without relying on (co)completeness or abelianness conditions. 

Section~2 is devoted to the fundamental tools and results in enriched category theory  that will be used throughout this work.  We fix a closed symmetric monoidal category $(\mathcal{V},\otimes,I,[-,-])$ with finite (co)limits. We denote by $\V_0$ its underlying ordinary category. Given a small $\mathcal{V}$-category $\mathcal{A}$ and a $\mathcal{V}$-category $\mathcal{B}$, we denote by $[\mathcal{A},\mathcal{B}]_0$ the ordinary category of $\mathcal{V}$-functors from $\mathcal{A}$ to $\mathcal{B}$ together with $\mathcal{V}$-natural transformations as morphisms. The convenience of working within the framework of enriched category theory is justified by Proposition~\ref{pop:free-V-cat}: Under suitable conditions, every small ordinary category $\mathcal{L}$ admits a $\V$-enrichment $\mathcal{L}_{\mathcal{V}}$, called \textit{the free $\mathcal{V}$-category}, such that, for every $\mathcal{V}$-category $\mathcal{B}$, there is an equivalence of categories
$
[\mathcal{L}_{\mathcal{V}},\mathcal{B}]_0
\cong
\operatorname{Fun}(\mathcal{L},\mathcal{B}_0),
$
where $\mathcal{B}_0$ denotes the underlying ordinary category of $\mathcal{B}$. In particular, for suitable choices of $\mathcal{L}$, we have the equivalence $\mathbf{C}_N(\mathcal{B}_0)\cong [\mathcal{L}_{\mathcal{V}},\mathcal{B}]_0$ (Example~\ref{ex:enriched-Qco}) and $\mathbf{C}_{\equiv_m}(\mathcal{B}_0)\cong [\mathcal{L}_{\mathcal{V}},\mathcal{B}]_0$ (Example~\ref{ex:m_periodic00}).

In Section~2, we continue recalling  relevant results on  monoid objects in $\V$ and left-objects over a monoid in $\V$. Our exposition draws extensively on the theory developed  in \cite{Par77,Par78}. We point out that  monoid objects and left-objects over a monoid may be regarded as general rings and left modules in $\V$, respectively; see Example~\ref{ex:unit_monoid_object}. We also note that the unit object $I$ carries a canonical commutative monoid structure, and there is a natural isomorphism ${}_I\mathcal{V}\cong\mathcal{V}$ (Example~\ref{ex:unit_monoid_object}). Another key notion in our study is that of dualizable $R$-left-objects (\eqref{internally_proj}), introduced in \cite{Bak20} in the context of  $R$-$R$-biobject. However, our study requires the corresponding notion for $R$-$S$-biobjects, where $R$ is a fixed,  and $S$ is an arbitrary monoid object. In Proposition~\ref{prop:dualiz_S_R}, we establish the equivalence between these two formulations. Dualizable $R$-left objects  may be regarded  as $\mathcal{V}$-analogues of finitely generated projective modules (Example~\ref{ex:dual_module}).

In Section~3, we restrict our study to the ordinary  $\V$-functor category $[\A, {}_R\V]$, where ${}_R\V$ denotes the $\V$-category of $R$-left objects  over a monoid object $R$ in $\V$ (\eqref{objects-over-monoid}).  
 In Proposition~\ref{prop:R-mod_valued}, we show that the structure of ${}_R\mathcal{V}$ can be encoded in the domain category $\mathcal{A}$. This allows us to focus our study on the $\mathcal{V}$-valued $\mathcal{V}$-enriched functor category $[\mathcal{A},\mathcal{V}]_0$.   We adapt the approach of   \cite{ECO26} to our setting. So  for each object $A\in\mathcal{A}$, we consider locally induced monoid object 
 $R_A:=\A(A,A)$ (Example~\ref{ex:local-monoid}), and the associated ordinary category 
  $({}_{R_A}\mathcal{V})_0$ of $R_A$-left-objects.  This approach has the advantage of covering  module categories over Morita context rings (Example~\ref{ex:morita_context}). In Proposition~\ref{prop:adj_q_ev_p}, we show that, for every object $A \in \A$, we have an adjoint triple
  $( { \rm q}_A, { \rm ev}_A, { \rm p}_A)$, which under suitable mild conditions (see \eqref{q_p_}), can be extended to an adjoint triple
 $$\xymatrix{  \prod_{A \in \A} ({}_{R_A} \V)_0 \ar@/^2em/[rrr]^{  { \rm q}} \ar@/_2em/[rrr]^{  { \rm p}} &&&    [\A, \V]_0 \ar[lll]_{ \ev} },
$$
where $\ev_A$ denotes the evaluation functor at $A$.

Starting from Subsection~\eqref{setup2}, we additionally assume that $\V_0$ is a weakly idempotent complete additive category.
We consider the exact structure $\mathcal{E}_{\textrm{lowt}} $ on $[\A, \V]_0$ given by: $\mathbb{E}$ belongs to $\mathcal{E}_{\textrm{lowt}} $ if and only if,  for every $A\in \A$, ${ \rm ev }_A(\mathbb{E})$ belongs to  the trivial  exact structure $\mathcal{E}_A^{\text{triv}}$ on $({}_{R_A} \V)_0$. This exact structure recovers the objectwise split exact structures on $\mathbf{C}_N(R)$ (Example~\ref{ex:R-enriched-N-chain}) and $\mathbf{C}_{\equiv_m}(R)$ (Example~\ref{ex:m_periodic22}). Under suitable mild conditions, $\mathcal{E}_{\textrm{lowt}} $-projective and $\mathcal{E}_{\textrm{lowt}} $-injective objects in $[\A, \V]_0$ are direct summands of $\V$-functors of the form
$$ { \rm q} ((M_A)_{A \in \A}), \quad \textrm{ and } \quad  { \rm p} ((M_A)_{A \in \A}), \quad  \quad (M_A)_{A \in \A} \in \prod_{A\in \A} ({}_{R_A} \V)_0,$$
respectively (Corollary~\ref{prop:proj_inj_trivial}). Moreover,  $([\A, \V]_0; \mathcal{E}_{\textrm{lowt}} )$ has enough projectives and injectives (Proposition~\ref{prop:ses1}). We let $\mathcal{P}$ and $\mathcal{I}$ denote the ideals of $\V$-natural transformations factoring through $\mathcal{E}_{\textrm{lowt}}$-projectives and $\mathcal{E}_{\textrm{lowt}}$-injectives, respectively. We thus obtain two stable categories 
$$\St_{ \mathcal{P}}([\A, \V]_0) \quad  \textrm{ and  }  \St_{ \mathcal{I}}([\A, \V]_0),$$
defined as the quotients of $[\mathcal{A},\mathcal{V}]_0$ by the ideals $\mathcal{P}$ and $\mathcal{I}$, respectively.
 These stable categories have a left and right triangulated category structures (see \eqref{Frob}), respectively. In Propositions~\ref{loop_suspension} and \eqref{loop_suspension2}, we provide explicit descriptions of the loop space $\Omega$ and the suspension space functor $\Sigma$.

Section 4 is devoted to proving our main results. We introduce two sets of  hypotheses on $\A$, Condition~\ref{conditions:frob} and Condition~\ref{conditions:frob2}, and prove the following: 

\vspace{3mm}

\textbf{Theorem~\ref{theo:proj-inj}:} 
\begin{enumerate}[i)]
\item If Condition~\ref{conditions:frob} is satisfied, then  $\mathcal{E}_{\textrm{lowt}} \mbox{-Proj} \subseteq  \mathcal{E}_{\textrm{lowt}} \mbox{-Inj}$.
\item If Condition~\ref{conditions:frob2} is satisfied, then  $ \mathcal{E}_{\textrm{lowt}} \mbox{-Inj} \subseteq  \mathcal{E}_{\textrm{lowt}} \mbox{-Proj}$.

\end{enumerate}

\vspace{3mm}

As a corollary of Theorem~\ref{theo:proj-inj} (Corollary~\ref{thm:frob}), if a $\V$-category $\A$ satisfies both  Conditions~\ref{conditions:frob} and \ref{conditions:frob2}, then $([\A, \V]_0;\mathcal{E}_{\textrm{lowt}})$ is a Frobenius exact category, so both stable categories coincide $\St_{ \mathcal{P}}([\A, \V]_0)=\St_{ \mathcal{I}}([\A, \V]_0)$ and carry a triangulated category structure.

Our second main result characterizes left and right homotopy to zero in $[\mathcal{A},\mathcal{V}]_0$ in terms of explicit families of $R_{\mathrm{F}(A)}$-morphisms. 
This description explains why, in the cases of $\mathbf{C}_N(R)$ and $\mathbf{C}_{\equiv_m}(R)$, the abstract notion of homotopy induced by the Frobenius exact structure coincides with the corresponding chain-homotopy relation, expressed in terms of families of morphisms. For a complete formulation, see Theorem~\ref{thm:homotopy}.

\vspace{4mm}

\textbf{Theorem~\ref{thm:homotopy}:} 
\begin{enumerate}[i)]
\item Suppose that $\A$ satisfies Condition~\ref{conditions:frob2}. A   $\V$-natural transformation in  
$[\A, \V]_0$  is right homotopic to zero if and only if it is   determined by a family $\{s_A\}_{A \in \A}$ of 
 $R_{{ \rm G}(A)}$-morphisms in $\V_0$ 
$$s_{A}: {}_{R_A}[R_{{ \rm G}(A)}, X(A)]  \longrightarrow Y({ \rm G}(A))$$  
satisfying certain equalities.
\item Suppose that $\A$ satisfies Condition~\ref{conditions:frob}. A $\V$-natural transformation  in $[\A, \V]_0$ is left homotopic to zero if and only if it is determined by  a family $\{s_A\}_{A \in \A}$ of 
 $R_{{ \rm F}(A)}$-morphisms in $\V_0$
$$s_{A}:X({ \rm F}(A)) \longrightarrow  R_{{ \rm F}(A)} \otimes_{R_A} Y(A)  $$  
satisfying certain equalities.
\end{enumerate}

\vspace{3mm}

Since the  formulation given in Theorem~\ref{thm:homotopy}  remains somewhat abstract in this general setting,   we consider the particular case in \eqref{particular_case_homotopy}, which is sufficiently general for our purposes; see Example~\ref{ex:chain_varying}.

Returning to Conditions~\ref{conditions:frob} and \ref{conditions:frob2}, one of the key hypotheses  is the existence of a right/left-weak Serre $\V$-functor on $\A$ (Definition~\ref{def:left_right_weak_Serre}). These functors play a crucial role in transition from the functor  ${ \rm p}_A$ to  ${ \rm q}_B$ for some $B \in \A$, and vice versa.

To motivate these  notions, let us consider the case $\mathcal{V}:=k\Mod$, where $k$ is a commutative ring and $\mathcal{A}$ is a small $k$-linear category. The formulation of right/left-weak Serre k-linear functors on $\A$ is similar to that of  the classical Serre functor. Recall that the classical Serre functor relative to a commutative ring $k$ on $\A$ (see \cite[Definition 4.1]{DSS17}) consists of a $k$-linear autoequivalence 
$F: \A \rightarrow \A$ together with a natural isomorphism between $k$-linear bifunctors
$$\sigma: \Hom_k(\Hom(-,F(?)),k) \longrightarrow \Hom(?,-),$$
where
$\A$ is  $\Hom$-finite, that is, $\Hom(A,B)$ is a finitely generated projective $k$-module for every objects $A$ and $B$ in $\A$. In contrast, for right/left-weak Serre k-linear functors, the locally induced endomorphism rings $R_A$ replace $k$, and the natural isomorphism in both variables is replaced by a natural retraction or section in a single variable. There are two main reasons for these modifications.   First, as already pointed out,  the $\mathcal{E}_{\textrm{lowt}}$-projective and $\mathcal{E}_{\textrm{lowt}}$-injective objects depend on the locally induced  endomorphism rings  \( R_A \), $A \in \A$.  If we considered $k$ in the natural transformations $\sigma_A$ and $\nu_A$  in Definition~\ref{def:left_right_weak_Serre}, we would need to impose additional conditions on the rings $R_A$ depending on $k$ to enable a transition from $k$ to $ R_A$. However, the locally induced rings in cases $\mathbf{C}_N(R)$ and $\mathbf{C}_{\equiv_m}(R)$ is $R$, where $R$ is an arbitrary ring. Second,
one interpretation of a classical Serre functor on a category $\A$ is that it provides a kind of `mirror reflection' of the category $\A$. So it forces having a certain degree of rigidity.  For example, even if $\A$ admits a classical Serre functor relative to $k$,  its $R$-extension $\A \otimes \mathcal{R}$ (\eqref{R-ext}) need not admit one for an arbitrary $k$-algebra $R$. For instance, let $\A$ be the category with a unique object whose endomorphism ring is $k$. Clearly, it admits a Serre functor relative to $k$. On the other hand, for a non-symmetric $k$-algebra $R$, its $R$-extension doesn't admit a Serre functor.  In contrast, right-weak and left-weak Serre $k$-linear functors are stable under such extensions; see Remark~\ref{remark:serre_R_extension}.

In Remark~\ref{Remark:classical_Serre}, we discuss two situations in which the existence of a classical Serre functor relative to $k$ on $\A$ implies the existence of right/left-weak Serre functors: when all endomorphism rings coincide with the base ring $k$, and when the Serre functor is the identity.
 By Remark~\ref{remark:serre_R_extension}, so does their $R$-extensions.  This provides a theoretical explanation for the cases of  $\mathbf{C}_N(R)$ (Example~\ref{ex:N_complex_Serre}) and $\mathbf{C}_{\equiv_m}(R)$ (Example~\ref{ex:m_periodic}). On the other hand, in Example~\ref{ex:non_Serre}, we provide a $k$-linear category $\A$ which admits both right/left-weak Serre functor, but no  classical Serre functor (Remark~\ref{remark:explan_non_Serre}). As explained in Remark~\ref{remark:explan_non_Serre}, it  does not arise as an $R$-extension of another well-behaved category. So the obstruction to the existence of a classical Serre functor lies solely in the requirement of having a natural isomorphism in both components. We note that this category satisfies both Conditions~\ref{conditions:frob} and \ref{conditions:frob2} (Example~\ref{ex:non_Serre2}). On the other hand,  in Examples~\ref{ex:right_repetitive_quiver} and \ref{ex:right-chain} with $R=k$, which admit right-weak Serre functors, although we have a natural isomorphism in both components,  the lack of a classical Serre functor stems from the requirement that $F$ be an autoequivalence.

Finally, in the Appendix, we address Question~\ref{question} proposed at the end of Section~4. in Proposition~\ref{prop:enriched_quotient}, we show that under the conditions given in Setup~\ref{setup3}, there exists a $\V$-category $\underline{[\A, \V]}$ whose underlying ordinary category is the projective stable category $\St_{ \mathcal{P}}([\A, \V]_0)$, that is,
$(\underline{[\A, \V]})_0 \cong\St_{ \mathcal{P}}([\A, \V]_0)$. As pointed out in Remark~\ref{remark_injective_stable}, replacing the hypothesis in Setup~\ref{setup3} with the assumption that  $\V_0$ admits $  \left| {\rm N}^+(A) \right |$-fold products for every $A\in \A$, we can construct a $\V$-category $\overline{[\A, \V]}$ whose underlying ordinary category is the injective stable category $\St_{ \mathcal{I}}([\A, \V]_0)$, that is,
$(\overline{[\A, \V]})_0 \cong\St_{ \mathcal{I}}([\A, \V]_0)$.

\section{Preliminares}
In this section, we review certain basic notions and results that  are essential for our study.

\subsection{Exact categories.} An {\em exact structure on an additive category $ \mathcal{C}$}   is an isomorphism-closed collection $\mathcal{E}$ of distinguished kernel-cokernel pairs  $(i,p)$ in $\mathcal{C}$
$$ \mathbb{E} : \xymatrix {  A \ar@{^{(}->}[r]^i & B \ar@{->>}[r]^p & C},$$
which are called \textit{admissible monic} and \textit{admissible epic} of $\mathbb{E}$, respectively,
satisfying:
\begin{enumerate}[(i)]
\item For every object $A \in \mathcal{C},$  the identity morphism $\id_A$ is both an admissible monic and an admissible epic;
\item Admissible monics and admissible epics are closed under compositions;
\item The class $\mathcal{E}$ is closed under  pushout  and pullback   along a morphism in $\mathcal{C}$.
\end{enumerate}

$\mathcal{C}$ equipped with an exact structure $\mcE$, denoted by $(\mathcal{C}; \mathcal{E})$, is called an  \emph{exact category}. A kernel-cokernel pair  $(i,p)$, or $\mathbb{E}$, in $\mcE$ is  called a {\em short exact sequence } in $(\mcC; \mathcal{E})$.  We refer to \cite{Buh} for a detailed treatment  on exact categories.


 \subsection{Quasi-abelian categories.}\cite[Definition~1.1.3]{Sch99}. An additive category $\mathcal{C}$ is called\textit{ quasi-abelian} if
 \begin{enumerate}[(i)]
\item Every morphism has a (co)kernel.
\item  The class of kernels is stable under pushout along arbitrary morphisms, and the class
of cokernels is stable under pullback along arbitrary morphisms.
 
 \end{enumerate}

 Of course, every abelian category is naturally quasi-abelian. However, the converse doesn't hold in general. The prototypical example of interest in this work is the category $\mathbf{Ban}$ of Banach spaces; see Example~\ref{ex:Banach}. The following two results are crucial for extending our results to quasi-abelian categories.
 
 \subsection{Proposition.}\label{prop:quasi_abelian}\cite[Proposition~1.1.7]{Sch99} (or \cite[Proposition~4.4]{Buh}). The class $\mathcal{E}_{\operatorname{max}}$ of all kernel-cokernel pairs in a  quasi-abelian category is an exact structure.

\vspace{5mm}

 A morphism $f:A \rightarrow B$ in an additive category with (co)kernels has the so-called \textit{analysis}
 \begin{equation}\label{eq:analysis_qa}
\xymatrix{   & A \ar[rrr]^-f \ar@{->>}[rd] &   &  & B \ar@{->>}[rd] & \\
 \Ker(f)\ar@{^{(}->}[ru] & & \operatorname{Coim}(f) \ar[r]^-{\widehat{f}} & \im(f)\ar@{^{(}->}[ru] & & \coker(f)},
  \end{equation}
where $\operatorname{Coim}(f):=\coker(\Ker(f)) $ and $\im(f):=\Ker(\coker(f))$.

\subsection{Proposition.}\label{prop:epi_monic_q.a}\cite[Proposition~4.8]{Buh}. Let $f$ be a morphism in a quasi-abelian category $\mathcal{C}$. The
canonical morphism $\widehat{f} : \operatorname{Coim}(f) \rightarrow \im(f) $ in the analysis of $f$  is both monic and epic.

\subsection{Relative projective and injective objects.} Let $(\mcC; \mcE)$ be an exact category. An object $P$ in $\mcC$ is said to be $\mathcal{E}$-\textit{projective} if for every short exact sequence  $\mathbb{E}$ in $\mcE$, the sequence $\Hom(P, \mathbb{E})$ of abelian groups is exact. The class of all $\mathcal{E}$-projective objects in $\mcC$ is denoted by $\mathcal{E}\mbox{-Proj}$. $(\mcC; \mcE)$ is said to have \textit{enough projectives }if, for every object $A$ in $\mcC$, there exists an admissible epic $P \twoheadrightarrow A$ in $\mathcal{E}$ with  $P \in \mathcal{E}\mbox{-Proj}$.

Dually, an $\mathcal{E}$-\textit{injective} object in  $\mcC$ is defined, and the class of $\mcE$-injective objects in $\mcC$ is denoted by $\mathcal{E}\mbox{-Inj}$. $(\mcC; \mcE)$ is said to have \textit{enough injectives} if for every object $A$ in $\mcC$ there exists an admissible monic $A \hookrightarrow E$ in $\mathcal{E}$ with  $E \in \mathcal{E}\mbox{-Inj}$.

\subsection{Frobenius categories.}\label{Frob}
Let $(\mcC; \mcE)$ be an exact category. Suppose that $(\mcC; \mcE)$ has enough projectives, and denote by $\mathcal{P}$ the ideal of morphisms in $\mcC$ which factors through an  $\mathcal{E}$-projective object. We let $\St_{\P}(\C)$ denote the stable category that has the same objects as $\mathcal{C}$ and morphisms modulo $\mathcal{P}$.  Suppose that 
to every object $A$ in $\mathcal{C}$ has been assigned a short exact sequence in $\mathcal{E}$
$$\xymatrix {  \Omega(A) \ar@{^{(}->}[r] & P(A) \ar@{->>}[r] & A},$$
where $P(A)$ is an $\mathcal{E}$-projective object. 
 By  \cite[Proposition~1.2, Theorem~2.12]{BM94}, we know that $\Omega$ induces an additive  functor 
\begin{align*}
\Omega: \quad & \St_{\P}(\C) \longrightarrow \St_{\P}(\C) ,
\end{align*}
which  is a part of a canonical left triangulated structure on $\St_{\P}(\C) $. $\Omega$ is called \textit{the loop space functor} of $\St_{\P}(\C) $. The left triangulated structure on $\St_{\P}(\C) $ is given by  the sequences $\xymatrix{\Omega(B) \ar[r] & \operatorname{cocone}(f) \ar[r] & A \ar[r]^-f & B}$, where the first two morphisms fit into the following pullback diagram
$$  
\xymatrix{   \Omega(B) \ar@{^{(}->}[r] \ar@{=}[d]& \operatorname{cocone}(f) \ar@{->>}[r] \ar[d] & A \ar[d]^f\\
\Omega(B) \ar@{^{(}->}[r] & P(B) \ar@{->>}[r] & B
}
$$
Two morphisms $f, g: A \rightarrow B$ in $\C$ are said to be \textit{left homotopic}, denoted by $f \sim^l g$, if the difference $f-g$ is trivial in $\St_{\P}(\C)$, or, equivalently, factors  through the morphism $P(B) \twoheadrightarrow B$.

Dually, if $(\mcC; \mcE)$ has enough injectives, and denote by $\mathcal{I}$ the ideal of morphisms in $\mcC$ which factors through an $\mathcal{E}$-injective object,  we let $\St_{\mathcal{I}}(\C)$ denote the stable category which has the same objects as $\mathcal{C}$ and morphisms modulo $\mathcal{I}$.  Suppose that 
to every object $A$ in $\mathcal{C}$ has been assigned a short exact sequence in $\mathcal{E}$
$$\xymatrix {  A \ar@{^{(}->}[r]& I(A) \ar@{->>}[r] &  \Sigma (A)}, $$
where $I(A)$ is an $\mathcal{E}$-injective object.  Similarly, $\Sigma$ induces an additive  functor 
\begin{align*}
\Sigma:  \quad & \St_{\mathcal{I}}(\C) \longrightarrow \St_{\mathcal{I}}(\C),
\end{align*}
which is a part of a canonical right triangulated structure. $\Sigma$ is called \textit{the suspension space functor} of $\St_{\mathcal{I}}(\C)$.  The right triangulated structure on $\St_{\mathcal{I}}(\C) $ is given by   the sequences
 $\xymatrix{A \ar[r]^-f & B \ar[r] & \operatorname{cone}(f) \ar[r] & \Sigma (A)}$, where the last two morphisms fit into the following pushout diagram
$$  
\xymatrix{   A \ar@{^{(}->}[r] \ar[d]_-f& I(A) \ar@{->>}[r] \ar[d] & \Sigma (A) \ar@{=}[d]\\
\Omega(B) \ar@{^{(}->}[r] & \operatorname{cone}(f) \ar[r] & \Sigma (A)
}
$$

 Two morphisms $f, g: A \rightarrow B$ in $\C$ are said to be \textit{right homotopic}, denoted by $f \sim^r g$, if the difference $f-g$ factors through an $\mathcal{E}$-injective object, or, equivalently, through the morphism $A \hookrightarrow I(A)$.


An exact  category  $(\mcC; \mcE)$ is said to be \textit{Frobenius}  if it has enough projectives and
injectives, and $\mathcal{E}\mbox{-Proj}=\mathcal{E}\mbox{-Inj} $. In this case, the stable categories $\St_{\mathcal{P}}(\C)$ and $\St_{\mathcal{I}}(\C)$ coincide, and the resulting category, denoted by $\St(\mathcal{C})$, is triangulated; see \cite[Section 13.4]{Buh}. Moreover, $\Sigma$ is an autoequivalence of $\St(\C)$, with inverse $\Omega$. Since the left and right homotopy relations on morphisms in $\mathcal{C}$ coincide,  we simply refer to them as \textit{homotopic}.

\section{Monoid objects  and objects over monoids.}
In this section, we summarize the basic notions and facts related to closed symmetric monoidal categories.

A \textit{monoidal category} consists of a category $\V$, a bifunctor
$
-\otimes- \colon \V \times \V \longrightarrow \V,
$
called the tensor product, a unit object $I \in \V$, and natural isomorphisms $a$ (associator), $l$ (left unitor), and $r$ (right unitor), satisfying the coherence axioms given in \cite[\S~1.1, Eqs.~(1.1), (1.2)]{Kelly}. A monoidal category $\V$ is called \textit{symmetric} if it is equipped with a natural isomorphism
$
c\colon A\otimes B \xrightarrow{\cong} B\otimes A,
$
called the \textit{ symmetry}, satisfying the additional coherence axioms. We write $(\V,\otimes,I)$ when referring to a (symmetric) monoidal category.

A symmetric monoidal category is said to be \emph{closed} if, for every $X \in \V$, the functor $-\otimes X \colon \V \rightarrow \V$ admits  a right adjoint $[X,-] \colon \V \rightarrow \V$.  In this  case,   $[-,-]$  defines  a bifunctor $\V^{\mathrm{op}} \times \V \rightarrow \V$, and it is called the \textit{internal hom} functor.  We denote a closed symmetric monoidal category by the quadruple $(\V,\otimes,I,[-,-])$.

We assume that the reader is familiar with the basic notions and results from enriched category theory, as outlined in \cite[Chapters.~1--2 (and parts of 3)]{Kelly}.  In particular, for a  symmetric monoidal category $(\V,\otimes,I)$, we will make frequent use of the notions of \emph{$\V$-categories}, \emph{$\V$-functors} between $\V$-categories, and \emph{$\V$-natural transformations}. Note that  there is a 2-functor $(-)_0: \V \mbox{-CAT} \longrightarrow \mbox{CAT}$
, called  the \textit{underlying ordinary category} functor; see \cite[\S~1.3]{Kelly}.

For a  $\V$-category $\A$ and two objects $A,B \in \A$, we denote the hom-object by $\A(A,B)$, and the composition rule on $\A$ by $\circ_{\A}$. $\A$ is said to be \textit{small} if its underlying category $\A_0$ is small in the usual sense. Given two $\V$-categories $\A$ and $\mathcal{B}$, with $\A$ small, the $\V$-functors from $\A$ to $\mathcal{B}$, together with the $\V$-natural transformations between them,  form an ordinary category, which we denote by $[\A, \mathcal{B}]_0$. Given two $V$-functors $X,Y : \A \rightarrow \V$, we denote the set of all $\V$-natural transformations from $X$ to $Y$  by $\V\mbox{-Nat}(X,Y)$.

As a note,  the internal hom functor on $\V$, if it exists, makes $\V$ naturally a $\V$-category; see \cite[\S~1.6]{Kelly}. 

The classical examples of  closed symmetric monoidal categories are the  categories  $(\mathbf{Set}, \times, \{ *\})$ and $(k\Mod, \otimes_k, k)$ of sets and $k$-modules over a commutative ring $k$, respectively. $\mathbf{Set}$- and $k\Mod$-enriched categories are precisely usual and $k$-lineal categories, respectively.  In the following, we
provide a list of   interesting (non)-abelian examples to which we will
apply the later results.



\subsection{Example.} \label{exa:Ch}
Let $k$ be a commutative ring, and let $\bfC(k)$ denote the category of chain complexes of $k$-modules. This category admits two distinct closed symmetric monoidal structures:
  
\begin{enumerate}[i)]
\item $(\bfC(k),\tTen,{ \rm S}^0(k),\tHom)$ where $\tTen$ denotes the \emph{total tensor product}, $\tHom$ the \emph{total Hom} of chain complexes, and ${ \rm S}^0(k) = 0 \to k \to 0$ is the $0$th \emph{sphere} complex with $k$ concentrated in degree $0$; see \cite[App.~A.2]{Chr}. A category enriched over this monoidal structure is known in the literature as a dg-category; see \cite[\S~2.2]{Keller}.

\item $(\bfC(k),\mTen,{ \rm D}^0(k),\mHom)$ where $\mTen$ denotes the \emph{modified} total tensor product, $\mHom$ the \emph{modified} total Hom, and ${ \rm D}^0(k) = 0 \to k \to k \to 0$ is the $0$th \emph{disc} complex concentrated in homological degrees $0$ and $-1$; see \cite[\S~2]{EG99}. 

\end{enumerate}

\subsection{Example.} 
 \label{exa:QcohX}
Let $X$ be any scheme.
  \begin{enumerate}[i)]
  \item $(\mathscr{O}_X \Mod,\otimes_{X},\mathscr{O}_X,\shHome_X)$ is a closed symmetric mo\-noi\-dal category, where $\mathscr{O}_X \Mod$ is the Grothendieck category of all sheaves of $\mathscr{O}_X$-mod\-ules on $X$;  see \cite[~II\S~5]{Hartshorne}. It is well-known that this is a Grothendieck category;  see \cite[Prop.~3.1.1]{AGr57}.

  \item \sloppy The category $\Qco(X)$ of quasi-coherent sheaves on $X$ is a Grothendieck  subcategory of $\mathscr{O}_X \Mod$; see 
  \cite[Cor.~(2.2.2)(i,ii)]{MR3075000} or 
  \cite[II Prop.~5.7]{Hartshorne}. 
  As $I = \mathscr{O}_X$ is quasi-coherent and quasi-coherent sheaves are closed under tensor products  (see \cite[Cor.~(2.2.2)(v)]{MR3075000}), it follows that $(\Qco(X),\otimes_{X},\mathscr{O}_X)$ is a monoidal subcategory of $(\mathscr{O}_X \Mod,\otimes_{X},\mathscr{O}_X)$. The restriction of the internal $\shHome_X$ on $\Qco(X)$ doesn't yield a quasi-coherent sheaf. However, using the so-called \emph{coherator} functor $Q_X $, which is the right adjoint of the inclusion functor \mbox{$\Qco(X) \hookrightarrow \mathscr{O}_X \Mod$}, the functor \smash{$\shHomqce_X := Q_X \circ \shHome_X$} yields a closed structure on $(\Qco(X),\otimes_{X},\mathscr{O}_X)$. 
  The category $\Qcoh{X}$ is Grothendieck by \cite[Lemma~1.3]{AJPV08}.
  
  Recall that if $X$ is an affine scheme, then $ {\Qco(X) \cong k \Mod}$ for some commutative ring $k$, and $ \otimes_{X} \cong \otimes_k$. 
  \end{enumerate}

\subsection{Example.}\label{ex:Banach}
Let $\mathbf{Ban}$ and $\mathbf{Ban}_1$ denote the categories of Banach spaces over the field $k$ ($\mathbb{R}$ or $\mathbb{C}$)  with bounded linear operators (equivalently, continuous linear operators) and with bounded linear operators of norm at most $1$, respectively. It is well known that $\mathbf{Ban}$ is a quasi-abelian category  with finite (co)limits  (see \cite[Chapter~IV, Remark~2.1.1]{Buh08t}), whereas $\mathbf{Ban}_1$ is a bicomplete non-additive category (see \cite[Chapter~IV, Proposition~1.4.1]{Buh08t}). Both categories $\mathbf{Ban}$ and $\mathbf{Ban}_1$ admit a closed symmetric monoidal structure
given by the \emph{projective tensor product}:  For two  Banach spaces $X$ and $Y$, consider their algebraic tensor
product $X\otimes_k Y$, that is, tensor product over $k$ as vector spaces. The \emph{projective norm} on
$X\otimes_k Y$ is defined by
\[
\|z\|_{\pi}
=
\inf\left\{
\sum_{i=1}^{n}\|x_i\|\,\|y_i\|
\;\middle|\;
z=\sum_{i=1}^{n}x_i\otimes y_i
\right\},
\]
where the infimum is taken over all finite representations of
$z\in X\otimes_kY$ as sums of elementary tensors. The
\emph{projective tensor product} $X\otimes_{\pi}Y$ is  defined as the
completion of $X\otimes_k Y$ with respect to the projective norm. The unit object of the projective tensor product is $k$. In both $\mathbf{Ban}_1$ and $\mathbf{Ban}$, the internal hom-object of $X$ and  $Y$ is the Banach space $\mathcal{L}(X,Y)$ of continuous linear operators from $X$ to $Y$, equipped with the operator norm
$$
\|T\|
=
\sup_{\|x\|=1}\|T(x)\|,
$$
see \cite[Example~6.1.9.h]{Borceux2}.

\subsection{Tensor product of $\V$-categories.}\label{tensor_product_V_cat} Let $\A$ and $\mathcal{B}$ be $\V$-categories. \textit{The tensor product of} $\A$ and $\mathcal{B}$, denoted by $\A \otimes \mathcal{B}$, is the $\V$-category with objects $\operatorname{ob}(\A) \times \operatorname{ob}(B)$ and with the hom-object
$$\A \otimes \mathcal{B}((A,B),(A',B')):=\A(A,A')\otimes \mathcal{B}(B,B')$$
together with the canonical composition and the identity element; see \cite[p. 12]{Kelly}

\subsection{} From now on, $\V$ denotes 
a fixed category with finite (co)limits, equipped with a fixed closed symmetric monoidal structure $(\otimes, I, [-,-])$. As already pointed out, the closed structure $[-,-]$ on $\V$ induces a natural $\V$-enrichment on $\V$ itself. To distinguish, we write $\V_0$ for the underlying ordinary category of $\V$, while keeping the notation $\V$ when it is viewed as $\V$-category.

\subsection{The enriched Yoneda lemma.}\label{Yoneda} 
Let $\A$ be a small $\V$-category $\A$, and $A$ be an object in $\A$. By \cite[\S~1.6, Eq.~1.29]{Kelly}, 
$\A(A,-):\A \rightarrow \V$ is a $\V$-functor, called the \textit{representable} $\V$-functor. For any $\V$-functor $X: \A \rightarrow \V$,  we always have the following $\V$-natural isomorphisms
$$
\V \mbox{-Nat}(\A(A,-), X) \cong \Hom_{\V_0}(I, X(A)),  \quad \quad \textrm{ \textit{the weak Yoneda lemma}.}
$$

Let  $\mathcal{B}$ be a $\V$-category. Given  a pair of  $\V$-functors $X,Y \colon \A  \rightarrow \mathcal{B}$, the \textit{end} of the $\V$-bifunctor $\mathcal{B}(X(-),Y(-))\colon \A \otimes \A^{\textrm{op}} \rightarrow \V$, if it exists, is an object  $K \in \V_0$  together with a family $ \lambda:=\{ \lambda_A\}_{A\in \A}$ of morphisms  $\lambda_A: K \rightarrow \mathcal{B}(X(A), Y(A))$ in $\V_0$ such that  for every pair of objects $A,A' \in \A$, the following  diagram commutes 
$$
\xymatrix{
K \ar[d]_-{\lambda_{A'}} \ar[rr]^-{\lambda_A} && \mathcal{B}(X(A), Y(A)) \ar[d]^{\rho_{AA'}}\\
\mathcal{B}(X(A'), Y(A')) \ar[rr]_-{\sigma_{AA'}}&& [\A(A,A'),\mathcal{B}(X(A), Y(A')) ]
}, 
$$
where the right vertical morphism $\rho_{AA'}$ and the bottom horizontal morphism $\sigma_{AA'}$ are induced, under the tensor-hom adjunction in $\V$, by the following two compositions of morphisms
$$\scalebox{0.88}{\xymatrix{\A(A,A') \otimes \mathcal{B}(X(A),Y(A)) \ar[r]^-{ Y \otimes  \id} & \mathcal{B}(Y(A),Y(A')) \otimes \mathcal{B}(X(A),Y(A)) \ar[r]^-{\circ_{\mathcal{B}}} & \mathcal{B}(X(A),Y(A')) }},$$
$$\scalebox{0.87}{\xymatrix{\A(A,A') \otimes \mathcal{B}(X(A'),Y(A')) \ar[r]^-{X\otimes \id} & \mathcal{B}(X(A),X(A')) \otimes \mathcal{B}(X(A'),Y(A')) \ar[r]^-{\circ_{\mathcal{B}}} & \mathcal{B}(X(A),Y(A')) }},$$
respectively, and $(K, \lambda)$ satisfies the universal property in the sense that if there exists another pair $(T, \beta)$ satisfying the same conditions, there exists a unique morphism $f: T \rightarrow K$ in $\V_0$ such that $\beta_A=\lambda_A \circ f$ for every $A \in \A$; see \cite[Eq.~1.40, \S~2.1]{Kelly}. In this case, the end  of $\mathcal{B}(X(-),Y(-))$ is denoted by 
$[\A, \mathcal{B}](X,Y)$\footnote{The standard notation for the end of a $\V$-bifunctor $T:\A \otimes \A^{\textrm{op}} \rightarrow \V$ is in fact $\int_{A\in \A} T(A, A)$.}. As noted in   \cite[\S~2.1, Eq.~2.2]{Kelly}, if $\V_0$ has $\vert \mbox{ob}(\A) \vert$-fold products\footnote{Products of any family of objects in $\V_0$  indexed by a set of cardinality at most  $ \vert \mbox{ob}(\A) \vert$ exist.}, then the end $[\A, \mathcal{B}](X,Y)$ exists, and is given by the equalizer of  
  the following two morphisms in $\V_0$ 
\begin{equation}\label{eq:end}
\xymatrix{\prod_{A \in \A} \mathcal{B}(X(A),Y(A)) \ar@<.5ex>[r] \ar@<-.5ex>[r]  & \prod_{A,A'\in \A}[ \A(A,A'), \mathcal{B}(X(A),Y(A'))]}.
\end{equation}
Here,  the two parallel morphisms are induced by  $\sigma_{AA'}$ and $\rho_{AA'}$, via the universal property of products. Of course, the condition of having $\vert \mbox{ob}(\A) \vert$-fold products is satisfied whenever $\V_0$ is complete.

 In the particular case $\mathcal{B}:=\V$, the end $[\A, \V](X,Y)$, if it exists,  is the object which for any object $T$ in $\V$, satisfies the following natural isomorphism
 \begin{equation}\label{eq:univ_end}
\Hom_{\V_0}(T, [\A, \V](X,Y)) \cong \V \mbox{-Nat}( X, [T,Y]  );
\end{equation}
where  $[T,Y]$ denotes the $\V$-functor $[T, -] \circ Y$; see \cite[Lemma~6.3.3]{Borceux2}. Equivalently, it is the object satisfying the  $\V$-natural isomorphism 
$$[T,[\A, \V](X,Y)] \cong [\A, \V](X, [T,Y]);$$
 see \cite[\S~2.2, Eq.~(2.13)]{Kelly}.

 The strong version of Yoneda lemma asserts that   for every $\V$-functor $X: \A \rightarrow \V$ and every object $A \in \A$, $[\A, \V](\A(A,-), X)$  always exists (without any additional assumptions on $\V_0$), and there is a $\V$-natural isomorphism
\begin{align*}
\
[\A, \V](\A(A,-), X) &\cong X(A), & \quad \textrm{ \textit{ the strong Yoneda lemma};}
\end{align*}
see \cite[\S~2.4]{Kelly}.  


Returning to the general case for $\mathcal{B}$, if for every pair of $\V$-functors $X,Y: \A \rightarrow \mathcal{B}$, the end $[\A, \mathcal{B}](X,Y)$ exists, then the  collection of all $\V$-functors from $\A$ to $\mathcal{B}$   can be endowed with a $\V$-category structure, denoted by 
$[\A, \mathcal{B}]$, given  by the hom-object $[\A,B](X,Y)$;
see \cite[\S~2.2]{Kelly}. Besides, the  underlying ordinary category $([\A, \mathcal{B}])_0$ of $\V$-category $[\A, \mathcal{B}]$ is precisely $[\A, \mathcal{B}]_0 $, that is, the morphisms in its underlying ordinary category   are precisley   $\V$-natural transformations
$$\Hom_{\V_0}(I, [\A,B](X,Y))=\V \operatorname{-Nat}(X,Y);$$
 see \cite[\S~2, Eqs.~2.12, 2.16]{Kelly}.

\vspace{5mm}

We recall the following simpler characterization of $\V$-functors from $\A$ to $\V$, obtained via the tensor-hom adjunction in $\V$. It will be used frequently in the upcoming sections.

\subsection{Proposition.}\label{X_AB}  Giving a $\V$-functor $X: \A \rightarrow \V$ is equivalent to giving a assignment from $\mbox{Ob}(\A)$ to $\mbox{Ob}(\V)$, denoted by $X(A)$, and for every pair of objects $A,B \in \A$, a morphism in $\V$
$$\overline{X}_{AB}: \A(A,B)\otimes X(A) \longrightarrow X(B),$$ satisfying 
$$ 
\xymatrix{
\A(B,C) \otimes \A(A,B) \otimes X(A) \ar[rr]^-{\circ_{\A} \otimes \id  } \ar[d]_{\id \otimes \overline{X}_{AB}  } && \A(A,C) \otimes X(A), \ar[d]^{\overline{X}_{AC}}\\
\A(B,C) \otimes X(B) \ar[rr]_{\overline{X}_{BC}} && X(C)
} \quad \quad \textrm{ \textit{the composition axiom}}
$$ 
$$
\xymatrix{I \otimes X(A) \ar[rr]^-{j_A\otimes \id} \ar[rrd]_{l} && \A(A,A) \otimes X(A),  \ar[d]^{\overline{X}_{AA} }\\
&& X(A)
} \quad \quad \textrm{ \textit{ the unit axiom}.}
$$ 
 It is also  equivalent to giving  the morphism  
$$\underline{X}_{AB}: X(A) \longrightarrow [\A(A,B),X(B)].$$
satisfying analogous  composition and  unit axioms. 
 
 \subsection{Proposition}\label{V-natural}\cite[Corollary 2.5]{HG16}.
 Let $X,Y: \A \rightarrow \V$ be two $\V$-functors. Giving a $\V$-natural transformation $\alpha: X \Rightarrow Y$ is equivalent to  giving a family of morphisms $\alpha_A: X(A) \rightarrow Y(A)$ in $\V_0$, $A \in \A$,  satisfying
$$
 \alpha_B \circ \overline{X}_{AB} =\overline{Y}_{AB} \circ (\alpha_A \otimes \id), \quad \quad A,B \in \A,
$$
or equivalently, 
$$
 \underline{Y}_{AB} \circ \alpha_A   =[\id, \alpha_B] \circ \underline{X}_{AB}, 
 \quad A,B \in \A.$$ 
  \vspace{5mm}
  
We recall the following result, which shows the convenience of working in the framework of enriched categories for obtaining results in functor categories with values in certain non-module categories ,  such as $\Qco(X), $ the category of quasi-coherent sheaves over a scheme $X$ (Example ~\ref{exa:QcohX}-(ii)) and $\mathbf{Ban}$, the category of Banach spaces (Example~\ref{ex:Banach}). 
\subsection{Proposition.}\label{pop:free-V-cat}\cite[\S~2.5]{Kelly}. 
Let $\mathcal{L}$ be a small ordinary category. If $\V_0$ has $  \vert \Hom(A,B) \vert$-fold coproducts for every objects $A$ and $B$ in $\A$, then  \textit{the free $\V$-category} $\mathcal{L}_\V$ of $\mathcal{L}$ defined with  the same objects as $\mathcal{L}$, and the hom-object 
$$\mathcal{L}_\V(A,B):=\bigoplus_{f \in \Hom_{\mathcal{L}}(A,B) }I_f, \quad \quad I_f:=I,$$
together with the canonical associative and unit morphisms is a $\V$-category. Moreover,  for every $\V$-category $\mathcal{B}$, 
 we have the following equivalence of categories
\begin{equation}\label{eq:functors_enriched_usual}
[\mathcal{L}_\V, \mathcal{B}]_0 \cong \mbox{Fun}(\mathcal{L}, \mathcal{B}_0),
\end{equation}
where $\mbox{Fun}(\mathcal{L}, \mathcal{B}_0)$ denotes the category of ordinary functors from $\mathcal{L}$ to $\mathcal{B}_0$.
\vspace{5mm}

\subsection{Remark.} We point out that the statement in \cite[\S~2.5]{Kelly} assumes  the (co)completeness of $\V_0$. However, as it is clear from its explanation, this assumption is required to construct $\V$-category $\mathcal{L}_\V$ for every small ordinary category $\mathcal{L}$, and consequently,  to obtain the left adjoint $(-)_\V$ of the underlying ordinary category 2-functor $(-)_0: \V \mbox{-CAT} \longrightarrow \mbox{CAT}$.   However, in this work, what we need is the equivalence of categories $[\mathcal{L}_\V, \mathcal{B}]_0 \cong \mbox{Fun}(\mathcal{L}, \mathcal{B}_0)$ for certain small ordinary categories $\mathcal{L}$, rather than for all small categories. The version given in  Proposition~\ref{pop:free-V-cat} allows us to study with categories such as $\mathbf{Ban}$, which  is  additive, and has  finite (co)limits.


\vspace{5mm}

The most-known example of the free $\V$-category construction  is the so-called \textit{preaditivization} of an ordinary category $\mathcal{L}$, which is precisely the free $\Ab$-category of $\mathcal{L}$. So if $\mathcal{B}$ is a preadditive category (equivalently, an $\Ab$-enriched category), then
$$ \mbox{Add}(\mathcal{L}_{\Ab}, \mathcal{B}) \cong \mbox{Fun}(\mathcal{L}, \mathcal{B});$$
see \cite[p. 11]{Mitchell}. 

\subsection{Example.}\label{ex:enriched-Qco}
Let $N\geq 2$ be an integer, and $\mathcal{L}$ be the category whose set of objects is $\Z$, and
$$\Hom(i,j):=
\begin{cases}
\{ \bullet \}, & i -(N-1) \leq j \leq i\\
\emptyset, & \textrm{ otherwise}
\end{cases}$$
By assumption, $\V_0$  has finite colimits, so we can apply Proposition~\ref{pop:free-V-cat}. The free $\V$-category $\mathcal{L}_{\V} $
has the same objects as $\mathcal{L}$,  and the hom-object between $i,j \in \Z$ is given by 
$$\mathcal{L }_{ \V} (i,j):=
\begin{cases}
I, & i -(N-1) \leq j \leq i\\
0, & \textrm{ otherwise}
\end{cases}$$
Moreover, for every $\V$-category $\mathcal{B}$, we  have the equivalence of categories
$$
[\mathcal{L}_{\V}, \mathcal{B}]_0  \cong \mbox{Fun}(\mathcal{L},\mathcal{B}_0)
$$  
which is also  equivalent to the category $\bfC_N(\mathcal{B}_0)$   of $N$-chain complexes of $\mathcal{B}_0$. Applying Example \ref{exa:QcohX},
we  have the equivalence of categories
\begin{align*}
[\mathcal{L}_{\mathscr{O}_X\operatorname{-Mod}}, \mathscr{O}_X\operatorname{-Mod}]_0 & \cong \mbox{Fun}(\mathcal{L},\mathscr{O}_X\operatorname{-Mod}) \cong \bfC_N(\mathscr{O}_X\operatorname{-Mod} )\\  
[\mathcal{L}_{\Qco(X)}, \Qco(X)]_0  & \cong \mbox{Fun}(\mathcal{L},\Qco(X)) \cong \bfC_N(\Qco(X)), 
\end{align*}
 Note that  if $X$ is an affine scheme, then  $\Qco(X)$ is equivalent to $k \Mod$ for some commutative ring $k$. In this case, the free $k\Mod$-category of $\mathcal{L}$ is denoted by $\mathcal{L }_{k}$. 

If we consider $\V:=\mathbf{Ban}$ (Example~\ref{ex:Banach}), then we obtain
$$[\mathcal{L}_{\mathbf{Ban}}, \mathbf{Ban}]_0  \cong \mbox{Fun}(\mathcal{L},\mathbf{Ban}) \cong \bfC_N(\mathbf{Ban}).$$
For the use of the usual chain complexes of Banach spaces, see \cite[\S~5.2]{Run02}.

\subsection{Example.}\label{ex:m_periodic00} Let $m$ be a positive integer, and  $\mathcal{L}$ be the category whose set of objects is
$\{0,1,\ldots, m-1\}$, and
$$\Hom(i,j):=
\begin{cases}
\{ \bullet \}, & 0 \leq i< m-1 \textrm{ and } j=i, i+1;\\
\{ \bullet \}, & i=m-1 \textrm{ and } j=m-1, 0;\\
\emptyset, & \textrm{ otherwise}
\end{cases}$$
Clearly, $\V_0$ satisfies the conditions  of Proposition~\ref{pop:free-V-cat} for $\mathcal{L}$. The free $\V$-category $\mathcal{L}_{\V} $
has the same objects as $\mathcal{L}$,  and the hom-object between $i,j$ is given by 
$$\Hom(i,j):=
\begin{cases}
I, & 0 \leq i< m-1 \textrm{ and } j=i, i+1;\\
I, & i=m-1 \textrm{ and } j=m-1, 0;\\
0, & \textrm{ otherwise}
\end{cases}$$
For every $\V$-category $\mathcal{B}$, we  have the equivalence of categories
$$
[\mathcal{L}_{\V}, \mathcal{B}]_0  \cong \mbox{Fun}(\mathcal{L},\mathcal{B}_0).$$

Note that when $\V:=\Ab$ and $\mathcal{B}$ is an abelian category, it is known that $\Add(\mathcal{L}_{\V},\mathcal{B})$ is equivalent to the category  $\bfC_{\equiv_m}(\mathcal{B})$ of $m$-periodic chain complexes of $\mathcal{B}$ (or \textit{$m$-cycle complexes}) as defined in \cite[\S~7]{PX97}). Motivated by this, for a more general  $\V$, we shall refer to the objects of $[\mathcal{L}_{\V}, \mathcal{B}]_0\equiv \bfC_{\equiv_m}(\mathcal{B}_0)$ as $m$-periodic chain complexes of $\mathcal{B}_0$.

\subsection{Monoid objects.}\cite{Par77}.
 An object $R$ in $\V$ is said to be \textit{monoid} if there exist two morphisms $ \mu_R: R \otimes R \rightarrow R$ and $\eta_R: I \rightarrow R$ in $\V_0$ satisfying the following commutative diagrams
$$\xymatrix{
R \otimes R \otimes R \ar[rr]^{\id \otimes \mu_R} \ar[d]_{\mu_R \otimes \id} && R \otimes R \ar[d]^{\mu_R}\\
R\otimes R \ar[rr]_{\mu_R} && R
}
\quad 
\xymatrix{I \otimes R \ar[rr]^{\eta_R \otimes  \id}  \ar[rrd]_r&& R \otimes R  \ar[d]{\mu_R}&& R \otimes I \ar[ll]_{\id \otimes \eta_R} \ar[dll]^l \\
 && R}$$
 $R$ is said to be \textit{commutative} if $\mu_R \circ c=\mu_R: R \otimes R \rightarrow R$. 
 
 The \textit{opposite monoid} of $R$, denoted by $R^\textrm{op}$, denotes  the triple $(R,\mu_{ R^\textrm{op}}, \eta_R )$, where $\mu_{ R^\textrm{op}}$ is the composition 
 $$\xymatrix{R\otimes R \ar[r]^{c} & R \otimes R \ar[r]^-{\mu_R} & R}.$$
 It is immediate that $R^\textrm{op}$ is again a monoid object in $\V$.

 Given two monoids $(R, \mu_R, \eta_R)$ and $(S, \mu_S, \eta_R)$  in $\V$, $R \otimes S$ is naturally a monoid object in $\V$ equipped with 
 $$\xymatrix{(R\otimes S) \otimes (R \otimes S) \ar[r]^-c & (R \otimes R) \otimes (S \otimes S) \ar[rr]^-{\mu_R \otimes \mu_S} && R\otimes S}  $$
 $$\xymatrix{ I \cong I \otimes I \ar[rr]^{\eta_R \otimes \eta_S} && R \otimes S} $$
A morphism $f: R \rightarrow S$ in $\V_0$ is said to be a \textit{monoid-morphism} if 
 $$f \circ \mu_R= \mu_S \circ (f \otimes f)\colon R \otimes R \rightarrow R \quad \textrm{ and }  \quad f \circ \eta_R= \eta_S \colon I \rightarrow S.$$
 
 \subsection{Objects over a monoid.}\label{objects-over-monoid} \cite{Par77}.
 Let $R$ be a monoid  in $\V$, an \textit{$R$-left-object} in $\V$ is an object $M \in \V$ together with a morphism $m: R \otimes M \rightarrow M$ in $\V_0$ satisfying the following commutative diagrams
 $$\xymatrix{
R \otimes R \otimes M \ar[rr]^{\mu_R \otimes \id} \ar[d]_{ \id \otimes m} && R \otimes M \ar[d]^m \\
R\otimes M \ar[rr]_{m} && M
}
\quad 
\xymatrix{I \otimes M \ar[rr]^{\eta_R \otimes  \id}  \ar[rrd]_r&& R \otimes M  \ar[d]^m\\
 && M}$$
 
Given two  $R$-left-objects $(M,m)$ and $(N,n)$ in $\V$, a morphism $f: M \rightarrow N$ in $\V_0$ is said to be a $R$-\textit{morphism}  if 
 $$ n \circ (\id \otimes f)=f \circ m \colon R \otimes M \longrightarrow N.$$
 The composition of $R$-morphisms in $\V_0$ is again an $R$-morphism. So 
 we get a usual category $({}_R \V)_0$  whose objects are the  $R$-left-objects in $\V$, and whose morphisms are the $R$-morphisms. We denote the set of $R$-morphisms from $M$ to $N$ by $\Hom_R(M,N)$.
 
The notions of  \textit{$R$-right-object} and of an \textit{$R$-morphism} between two $R$-right-objects in $\V$ are defined in a dual manner. Note that an $R$-right-object is precisely an  $R^\textrm{op}$-left-object.  Similarly, given two  monoid objects $R$ and $S$ in $\V$, $R$-$S$-\textit{biobjects} in $\V$ and \textit{$R$-$S$-morphisms} are defined, and they are precisely $(R\otimes S^{\textrm{op}})$-left objects. We denote by  $(\V_R)_0$ and $({}_R \V_S)_0$ the usual categories of  $R$-right-objects  and $R$-$S$-biobjects in $\V$, respectively.

\subsection{Example.}\label{ex:unit_monoid_object} The unit object $I$ of $\V$ with the left unitor $l_I$  is a commutative  monoid object in $\V$, and there exists the equivalence of categories $({}_I \V)_0 \cong \V_0 \cong (\V_I)_0$; see \cite[Corollary 2.3]{Par77}.   
 

 \subsection{Example.}Let $k$ be a commutative ring. A  $k$-module $R$  is a (commutative) monoid object in $(k \Mod, \otimes_k, k)$ if and only if it is a (commutaive) $k$-algebra; see \cite[p.~196]{Par77}.  In  particular, when $k:=\Z$, the (commutative) monoid objects in $\Ab$ are precisely the usual (commutative) rings, while  $R$-left-objects in $\Ab$ are precisely the usual  left $R$-modules.

 \subsection{Example.} An associative unital Banach algebra $A$ with  $\| 1_A \|=1$ is precisely a monoid object in the closed symmetric monoidal category $\mathbf{Ban}_1$; see \cite[p.~197]{Par77}. On the other hand, a non-trivial monoid object $A$ in $\mathbf{Ban}$ is a Banach space, which is a $k$-algebra,  such that there exist  positive real numbers $M $ and $C$ such that for every $a,b \in A$, $\| a \cdot b \| \leq M \| a \| \| b\|$   and $\| 1_A \| =C$.
\subsection{Example.}\label{ex:local-monoid} Let $\A$ be a $\V$-category.  It is straightforward to verify that, for any object $A$ in $\A$, the hom-object $\A(A,A)$ together with the canonical morphisms induced from the $\V$-structure of $\A$
$$\circ_{\A}: \A(A,A) \otimes \A(A, \A) \longrightarrow \A(A,A), \quad j_A:I \longrightarrow \A(A,A)$$
is a monoid object in $\V$, which will be denoted by $R_A$. By the composition rule in $\A$, for any two objects $A$ and $B$ in $\A$, $\A(A,B)$ is  an $R_B$-$R_A$-biobject in $\V$.

\subsection{Proposition.}\label{prop:V_R}\cite[Theorem 2.2, Corollary 2.4]{Par77}. 
 For any monoid $R $ in $\V$, the forgetful functor $U: ({}_R \V)_0 \longrightarrow \V_0$ has the left adjoint functor $R \otimes -$, and creates (co)limits. Besides, if $\{G_i \}_{i \in I}$ is a generating set in $\V_0$, so is $\{R \otimes G_i \}_{i \in I} $ in $({}_R \V)_0$. 
 
 \vspace{5mm}
 
 The following corollary follows  directly from Proposition~\ref{prop:V_R}.
\subsection{Corollary.}\label{cor:V_R_quasi-abelian}  For any monoid $R $ in $\V$, $({}_R \V)_0$ has finite (co)limits. Additionally, if $\V_0$ is (quasi-)abelian (or Grothendieck), so is $({}_R \V )_0$.

\subsection{$\V$-structure on $R$-$S$-biobjects.}\label{V_cat_R}\cite[\S~3]{Par77}. Let $R$, $S$ and $T$ be  monoids in $\V$. Given an $R$-$S$-biobject $(M,m) $ and an $S$-$T$-biobject $(N,n)$   in $\V$, respectively, we consider the following pair of  morphisms of $R$-$T$-biobjects
$$\xymatrix{M \otimes S \otimes N \ar@<1ex>[rr]^{m \otimes \id} \ar@<-1ex>[rr]_{\id \otimes n} && M\otimes N \ar@{.>>}[r] & M \otimes_S N}.$$
Since $({}_R\V_S)_0$ has finite colimits, their coequalizer exists, and is denoted by $M\otimes_S N$. 

By a dual construction, given an $R$-$S$-biobject $(M,m) $ and an $R$-$T$-biobject $(N,n)$   in $\V$, the object ${}_R[M,X]$ is defined as the equalizer of two  $S$-$T$-morphisms from $[M,N]$ to $[R \otimes M, N]$, which, under the tensor-hom adjunction in $\V$, correspond to
$$\xymatrix{R \otimes M \otimes[M,N] \ar[rr]^-{\id \otimes \xi_N}&& R \otimes N \ar[r]^-n & N},$$ 
$$\xymatrix{R \otimes M \otimes[M,N] \ar[rr]^-{m \otimes \id}&& M \otimes [M,N] \ar[r]^-{\xi_N} & N},$$  where $\xi_N$ denotes the counit morphism of  the tensor-hom adjuntion.

By the universal property of (co)equalizers, these constructions induce the following  functors\footnote{They are additive whenever $\V_0$ is an additive category.}
 \begin{align}
 - \otimes_S -:& ({}_R \V_S)_0 \times ({}_S \V_T)_0 \longrightarrow ({}_R \V_T)_0 \label{functors1}\\
 {}_R[-,-]: & ({}_R \V_S)_0^{\textrm{op}} \times ({}_R \V_T)_0  \longrightarrow ({}_S \V_T)_0\label{functors1*} 
   \end{align}
  satisfying the following natural isomorphisms
\begin{align}
{}_R [M\otimes_S N,X] \cong\ & {}_S [N, {}_R[M,X] ].\label{adj_functors2}\\
  \Hom_R(M\otimes_S N,X) \cong & \Hom_S(N, {}_R[M,X]  ) \label{adj_functors2*} 
 \end{align}
for every $R$-$S$-biobject $M$, $S$-$T$-biobject $N$, and $R$-left-object $X$; see \cite[Propositions~3.10-3-.11]{Par77}. Note that the isomorphism \eqref{adj_functors2} is in fact an isomorphism as $T$-left-objects. 

 By  \cite[Corollary 3.3]{Par77}, the functor ${}_R[-,-]$  endows the class of $R$-left-objects in $\V$ with a $\V$-category structure, which we denote  by ${}_R \V$. Its underlying ordinary  category is precisely $({}_R \V)_0$.  Furthermore, the functors in \eqref{functors1} and \eqref{functors1*} are $\V$-functors, as well.

Since every right-object over a monoid is a left-object over its opposite monoid, the corresponding right-handed version $[-,-]_R$ of ${}_R[-,-]$. and the previous statements and  isomorphisms about right-objects over a monoid follow immediately.



\subsection{Dualizable objects in ${}_R \V$.}\label{internally_proj}
Recall from \cite[(6.5)]{HO23} that an object $M$ in $\mathcal{V}$ is called \textit{dualizable object} if, for every  $N,C \in \V_0$, the canonical morphism 
$$[M,N]\otimes C \rightarrow [M,N \otimes C]$$
 in $\V_0$, induced from tensor-hom adjunction,   is an isomorphism. In \cite{Bak20}, the author introduces a noncommutative version of this notion for  left/right-objects over a monoid. Let $R$ and $S$ be  two   monoid objects  in $\V$. Suppose that $M$, $N$ and $C$ are, respectively,  an $R$-left-object, an $R$-$S$-biobject, and an $S$-left-object in $\V$. Applying the adjunction $(M\otimes-,{}_R[M,-])$ between $\V_0$ and $({}_R\V)_0$ to  the canonical $R$-morphism $\xi_N \otimes \id: M \otimes ({}_R[M, N]) \otimes_S C \rightarrow N \otimes_S C,$
where $\xi_N$ is the counit\footnote{As already pointed out in \eqref{adj_functors2*},  the pair $(M \otimes-,{}_R[M,-])$ is also an adjoint pair between $(\V_S)_0$ and $({}_R \V_S)_0$.} morphism of the adjunction, 
we obtain a $\V$-natural morphism in $\V_0$
\begin{equation}\label{evalut_monoid}
{}_R[M,N] \otimes_S C \longrightarrow  {}_R[M,N \otimes_S C].
\end{equation}

Similarly, given an $R$-left-object M,   an $R$-$S$-biobject $N$, and  an $S$-right-object $C$, there is a $\V$-natural morphism in $\V_0$
\begin{equation}\label{evaluat_monoid2}
[N,C]_S \otimes_R M   \longrightarrow [ {}_R[M,N],C]_S.
\end{equation}

In \cite[Definition 7, Lemma 5-(8)]{Bak20}, the author  considers the case $R=S$, and calls $M$ \textit{dualizable in ${}_R \V$} if the morphism \eqref{evalut_monoid} is an isomorphism  in $\V_0$ for every $R$-$R$-biobject $N$ and every $R$-left-object $C$. By Example~\ref{ex:unit_monoid_object}, we already know that  the unit object $I$ is a commutative monoid in $\V$, and every object in $\V$ is naturally an $I$-$I$-biobject. So the classical notion of dualizable object in $\V$ coincide with the notion  of dualizable object in ${}_I\V$ in the sense of \cite{Bak20}.  However, to obtain the main results of Section~\ref{section:homotopy} (Theorem~\ref{theo:proj-inj} and Theorem~\ref{thm:homotopy}), we need a version of this notion that allows the morphism \eqref{evalut_monoid} to be an isomorphism for possibly different monoids $R$ and $S$. For it, we prove the following result.

\subsection{Proposition.}\label{prop:dualiz_S_R}
An $R$-left-object $M$ is dualizable in ${}_R \V$ if and only if, for every monoid $S$, every $R$-$S$-biobject $N$ and every $R$-left-object $C$, the morphism  \eqref{evalut_monoid} is an isomorphism.

\begin{proof}
The sufficiency part is immediate. Now we prove the necessity part. Suppose that $M$ is dualizable in ${}_R \V$. Taking $M^*:={}_R[M,R]$, which is a $R$-right-object,  we have the $\V$-natural isomorphism of $\V$-functors
$$ M^* \otimes_R- \Longrightarrow {}_R[M,R \otimes_R-] \cong {}_R[M,-]: {}_R \V \longrightarrow \V.$$

Recall that  for a given  $R$-$S$-biobject $N$, the $S$-right-object structures on  both $M^*\otimes_R N$ and ${}_R[M,N]$ are given precisely by
$$ \xymatrix{(M^*\otimes_R N) \otimes S  \ar[r]^-{\cong} & M^*\otimes_R (N \otimes S) \ar[r]^-{\id \otimes n} & M^*\otimes_R N }$$
$$ \xymatrix{{}_R[M,N]\otimes S \ar[r] & {}_R[M,N \otimes S] \ar[r]^-{[\id, n]} & {}_R[M,N] }.$$
So the isomorphism $M^* \otimes_R N \longrightarrow  {}_R[M,N]$ is in fact an $S$-morphism, and therefore an isomorphism of  $S$-right objects.  Thus,  for any $S$-left object $C$, we have the following commutative diagram
$$
\xymatrix{
(M^* \otimes_R N) \otimes_S C \ar[r]^{\cong} \ar[d]_{\cong} & M^* \otimes_R(N \otimes_S C) \ar[d]^{\cong}\\
 ({}_R[M,N]) \otimes_S C \ar[r] & {}_R[M,N \otimes_S C]
}
$$ 
The top horizontal arrow is an isomorphism by associativity, while the two vertical arrows are isomorphisms by the preceding observation and the dualizability of $M$. Hence, the bottom horizontal arrow is also an isomorphism.
\end{proof}

\subsection{Proposition.}\label{prop:dualizable_codualizable}
Let $M$ be an $R$-left object. If $M$ is a dualizable object in ${}_R\V$, then for every monoid $S$, every $R$-$S$-biobject $N$ and every $S$-right-object $C$, the morphism \eqref{evaluat_monoid2} is an isomorphism.

\begin{proof}
By \cite[Lemma 6]{Bak20}, $M^*={}_R[M,R]$ is a dualizable object in $\V_R$, and $M \cong M^{**}=(M ^*)^*$.  So for an  $R$-$S$-biobject $N$ and an $S$-right-object $C$, we have
\begin{align*}
 [N,C]_S \otimes_R M  \cong [N,C]_S \otimes_R (M ^*)^* & \cong [M^*, [N,C]_S ]_R\\
 & \cong [M^* \otimes_R N, C]_S\\
 & \cong [  {}_R[M, N], C]_S
\end{align*}

\end{proof}

Since ${}_R[R,-] \cong \id\cong R\otimes_R-: {}_R \V \longrightarrow {}_R \V$ (see \cite[bottom of p. 183]{Par78}), any $R$-left-object that is a  direct summand of $R^n$, $n \geq 1$, is a dualizable object in ${}_R \V$. Such an $R$-left-object can be viewed as  `\textit{internally projective}'. 
 In the following result, we assume  that $\V_0$ is  quasi-abelian, and hence, by Corollary~\ref{cor:V_R_quasi-abelian}, so is $({}_R \V)_0$. Under this condition, 
 by the exactness of a functor, we refer to the exactness with respect to the exact structure $\mathcal{E}_{\operatorname{max}}$ given in Proposition~\ref{prop:quasi_abelian} on each corresponding category.

\subsection{Proposition.}\label{prop:internally_projective}
Let $M$ be an $R$-left-object.  Suppose that $\V_0$ is quasi-abelian. If $M$ is a dualizable object in ${}_R \V$, then it is internally projective, 
 that is, the functor 
$${}_R[M,-]: ({}_R \V)_0 \longrightarrow \V_0$$
 is exact.
 
\begin{proof}
It follows from the $\V$-natural isomorphism ${}_R[M,-]\cong M^*\otimes_R-$, which preserves (co)limits.
\end{proof}

\subsection{Remark.}
 The notion of internally projectivity  is introduced  in \cite{Lew99}, and  is formulated   for objects in $\V$, where $\V_0$ is assumed to be abelian. The aforementioned version is a natural generalization of it to the noncommutative quasi-abelian setting.

\subsection{Example.}\label{ex:dual_module} Let $k$ be a commutative ring. A $k$-module $M$ is a dualizable object in $k \Mod$ if and only if it is a finitely generated projective $k$-module; see \cite[p. 121]{LewisMay}.

Let $R$ be a $k$-algebra. A left $R$-module $M$ is  dualizable   if and only if  it is a finitely presented  projective left $R$-module. The necessity follows from Proposition~\ref{prop:internally_projective}. For the sufficiency, see \cite[p. xxxiii]{Kra21}. 


\subsection{Example.} Consider Example~\ref{exa:QcohX}. A quasi-coherent sheaf $\mathcal{F}$ over a scheme  $X$ is a dualizable object in $\Qco(X)$ if and only if it is  a  locally free sheaf of finite rank; see \cite[Proposition~4.7.5]{Bra14}.

\subsection{Example.} A Banach space $X$  is a dualizable object  in $\mathbf{Ban}$ if and only if it is finite-dimensional as a vector space. Indeed, in \cite[p. 4426]{ST12}, the category $\mathrm{TV}$ of Hausdorff, complete locally convex topological vector spaces equipped with the projective tensor product is considered. It is noted that an object $X$ in $\mathrm{TV}$ is dualizable if and only if it is finite-dimensional. Since $\mathbf{Ban}$ is a full monoidal subcategory of $\mathrm{TV}$, dualizability of a Banach space $X$ in $\mathbf{Ban}$ is equivalent to its dualizability in $\mathrm{TV}$.

\section{The left and right triangulated structures on Functor categories}

Throughout, 
we let $\A$ denote  a small $\V$-enriched category. We continue to  assume that $\V_0$ has finite (co)limits. In this section, we  focus on  the category $[\A, \V]_0$ equipped with the object-wise trivial exact structure.   In certain results, to obtain analogous statements for its enriched counterpart $[\A,\V]$, we would need to assume that  $[\A, \V]$ exists, that is, that the $\V$-structure on $\V$-functors from $\A$ to $\V$ exists; see \eqref{Yoneda}.

Before going on  further, we begin by proving Proposition \ref{prop:R-mod_valued}, which shows that if ${}_R \V$ appears as the codomain of a $\V$-functor category, then  its structure can be encoded in   the domain category, and therefore, we may  restrict our study to $[\A, \V]$.

\subsection{$R$-extension of a $\V$-category.} \label{R-ext}
Let $R$ be a monoid object in $\V$. We let $\mathcal{R}$ denote the $\V$-category with a single object $0$, whose  hom-object is $\mathcal{R}(0,0):=R$ with the composition rule and the identity morphism given by  $\mu_R:R \otimes R \rightarrow R$ and $\eta_R: I \rightarrow R$, respectively. Recall from \eqref{tensor_product_V_cat} the tensor product $\V$-category $  \A \otimes \mathcal{R}$. Its objects are pairs $(A,0)$, where $A \in \operatorname{ob}(\A)$, and its hom-objects are given by
$$( \A \otimes \mathcal{R})((A,0),(B,0)):= \A(A,B) \otimes R, \quad A,B \in \A.$$
For simplicity, we shall continue to denote the objects $(A,0)$ of $\A\otimes\mathcal{R}$ simply by $A$. So we have
$$( \A \otimes \mathcal{R})(A,B):= \A(A,B) \otimes R, \quad A,B \in \A.$$
The composition and identity morphisms are given, respectively, by
$$\scalebox{0.85}{\xymatrix{(  \A(B,C)\otimes R) \otimes ( \A(A,B)\otimes R ) \ar[r]^-c &  ( \A(B,C) \otimes \A(A,B))\otimes (R \otimes R) \ar[rr]^-{ \circ_{\A}\otimes \mu_R } &&  \A (A,C)\otimes R }}$$
$$ \xymatrix{ I \ar[r]^-{\cong} & I \otimes I \ar[rr]^-{ j_A \otimes \eta_R } &&  \A(A,A) \otimes R}.$$

\subsection{Proposition.}\label{prop:R-mod_valued} For every monoid object $R$ in $\V$,  we have the equivalence of categories
$$[\A, {}_R \V]_0 \cong [ \A \otimes \mathcal{R}, \V]_0.$$
If in addition, $[\A, \V]$ exists, then  we have the equivalence of  $\V$-functors
$$[\A, {}_R \V] \cong [\A \otimes \mathcal{R}, \V].$$
\begin{proof}
Since the $\V$-category $\mathcal{R}$ associated to $R$ has only one object,  the $\V$-category $[\mathcal{R}, \V]$ exists. By \cite[\S~2.3]{Kelly} (see also \cite[Eqs.~(2.16) and (2.20)]{Kelly}), we have the equivalence of categories
$$[\A \otimes \mathcal{R}, \V]_0 \cong [\A, [\mathcal{R}, \V]]_0.$$ 
It remains to show the $\V$-equivalence $[\mathcal{R}, \V] \cong {}_R\V$. This is straightforward to verify (see also \eqref{evaluation}).

If $[\A, \V]$ exists, then so does $[\A, {}_R \V]$, since the forgetful functor $U: ({}_R \V)_0 \rightarrow \V_0$ creates (co)limits. By \cite[\S~2.3, Eq.~(2.29)]{Kelly}, we have the equivalence of $\V$-functors
$$[\A, {}_R \V] \cong [\A \otimes \mathcal{R}, \V].$$

\end{proof}

\subsection{Example.}\label{ex_N_chain_N}  Let $R$ be a monoid object in $\V$. Consider the free $\V$-category $\mathcal{L}_{\V}$ given in  Example~\ref{ex:enriched-Qco}. We already have the equivalence of  categories $[\mathcal{L}_{\V}, {}_R \V ]_0 \cong \operatorname{Fun}(\mathcal{L}, ({}_R \V)_0) \cong \mathbf{C}_N(({}_R \V)_0)$. By Proposition~\ref{prop:R-mod_valued}, it is also equivalent to the category $[ \mathcal{L}_{\V} \otimes \mathcal{R}, \V]$, where the hom-objects  in $\mathcal{L}_{\V} \otimes \mathcal{R}$ are given by 
$$(\mathcal{L}_{\V} \otimes \mathcal{R})(i,j):=
\begin{cases}
R, & i -(N-1) \leq j \leq i\\
0, & \textrm{ otherwise}
\end{cases}$$

\subsection{Example.}\label{ex_m_periodic11}  Let $R$ be a monoid object in $\V$. Consider the free $\V$-category $\mathcal{L}_{\V}$ given in  Example~\ref{ex:m_periodic00}. We have the equivalence of  categories 
$$[\mathcal{L}_{\V}, {}_R \V ]_0 \cong \operatorname{Fun}(\mathcal{L}, ({}_R \V)_0) \cong \mathbf{C}_{\equiv_m}(({}_R \V)_0) \cong [ \mathcal{L}_{\V} \otimes \mathcal{R}, \V],$$  where the hom-objects  in $\mathcal{L}_{\V} \otimes \mathcal{R}$ are given by 
$$(\mathcal{L}_{\V} \otimes \mathcal{R})(i,j):=
\begin{cases}
R, & 0 \leq i< m-1 \textrm{ and } j=i, i+1;\\
R, & i=m-1 \textrm{ and } j=m-1, 0;\\
0, & \textrm{ otherwise}
\end{cases}$$

\subsection{}\label{evaluation} In view of  Proposition~\ref{prop:R-mod_valued}, from now on we shall focus on the $\V$-valued $\V$-enriched functors on $\A$.

As already pointed out in Example~\ref{ex:local-monoid}, for any object $A$  in $\A$, the hom-object $R_A:=\A[A,A]$ is naturally a  monoid object in $\V$, and  $\A(A,B)$ is an $R_B$-$R_A$-biobject in $\V$. By  Proposition~\ref{prop:V_R} and \eqref{V_cat_R}, the collection of $R_A$-left-objects in $\V $ is a $\V$-enriched category  with hom-object ${}_{R_A}[-,-]$, whose underlying ordinary category $ ({}_{R_A}\V)_0$ has finite (co)limits.

 Let  $X: \A \rightarrow \V$ be a $\V$-functor. Recall from Proposition~\ref{X_AB} that  for every pair of objects $A,B \in \A$, there exist morphisms in $\V_0$
$$\overline{X}_{AB}: \A(A,B)\otimes X(A) \longrightarrow X(B),$$ 
$$\underline{X}_{AB}: X(A) \longrightarrow [\A(A,B),X(B)].$$
satisfying the  composition and  unit axioms.  One consequence of the previous data is that when $A=B$, the morphism $\overline{X}_{AA}$ endows $X(A)$ with the structure of a $R_A$-left-object, and  $\overline{X}_{AB}$ and $\underline{X}_{AB}$ are $R_B$- and $R_A$-morphisms, respectively. Moreover, by the universal property of (co)equalizers, there exist unique $R_B$- and $R_A$-morphisms, respectively,
\begin{align}
\overline{X}_{AB}: & \A(A,B) \otimes_{R_A} X(A) \longrightarrow X(B), \label{morph_X1}\\ 
\underline{X}_{AB}: & X(A) \longrightarrow {}_{R_B}[\A(A,B),X(B)],\label{morph_X2}
\end{align}
see \eqref{V_cat_R} for the construction of $\otimes_{R_A}$ and ${}_{R_A}[-,-]$.

 On the other hand, for a given object $A$ in $\A$, we consider the full $\V$-subcategory $\mathcal{R}_A \hookrightarrow \A$ (see \eqref{R-ext}) given by the monoid object $R_A$. This gives rise to the evaluation functor at $A$
 $$ \xymatrix{[\A, V]_0 \ar[rr]^-{\ev_A} && [\mathcal{R}_A, \V]_0 \cong ({}_{R_A}\V})_0.$$
It is worth emphasizing that if $[A, \V]$ exists, then  the evaluation functor at $A$  is in fact a $\V$-functor  
$$ \xymatrix{[\A, V] \ar[rr]^-{\ev_A} && [\mathcal{R}_A, \V] \cong {}_{R_A}\V}.$$ 


Using the $\V$-functors  in \eqref{functors1} and \eqref{adj_functors2*}, for a given $R_A$-left object $M$ in $\V$, we define the following $\V$-functors from $\A$ to $\V$ by 
\begin{align*}
{ \rm q}_A(M)(-) & := \A (A,-) \otimes_{R_A} M, \quad (-\otimes_{R_A}M) \circ \A (A,-)\\
{ \rm p}_A(M)(-) &:={}_{R_A}[ \A(-,A),M], \quad {}_{R_A}[ -,M] \circ \A(-,A).
\end{align*}
They induce two functors
$${ \rm q}_A,{ \rm p}_A: ({}_{R_A} \V)_0 \longrightarrow [\A, \V]_0.$$
Moreover, if $[\A, \V]$ exists, then they are $\V$-functors, as well.
Indeed, for a given object $A$ in $\A$, the $\V$-functor $\A(A,-): \A \rightarrow \V$ has image in $\V_{R_A}$, and therefore, under the $\V$-isomorphism $[\A \otimes {}_{R_A} \V, \V] \cong [ {}_{R_A} \V , [\A, \V]]$, the following composition of $\V$-functors
$$\xymatrix{\A \otimes {}_{R_A} \V \ar[rr]^{\A(A,-) \otimes \id} && \V_{R_A} \otimes {}_{R_A} \V \ar[rr]^-{- \otimes_{R_A} -} && \V} $$
induces the $\V$-functor ${ \rm q}_A  \colon {}_{R_A} \V \longrightarrow [\A, \V]$. Similar arguments can be applied to show that ${ \rm p}_A$ is a $\V$-functor.

The following result can be seen as an $R_A$-left-object analogue of the weak and strong Yoneda lemmas. By the adjointness in the enriched sense, we mean  the $\V$-natural isomorphism
 $$[\A,\V]({\rm q}_A(M),X) \quad \textrm{ and } \quad [\A,\V](X,{\rm p}_A(M)).$$

\subsection{Proposition.}\label{prop:adj_q_ev_p} For each object $A \in \A$, the  triple $({ \rm q}_A,\ev_A,    { \rm p}_A)$ is an adjoint triple in both  the ordinary and $\V$-enriched  senses. 

\begin{proof}
We only show the adjointness of the pair $({\rm q}_A,\ev_A )$, since the other one can be proven similarly. 

First, we show that the pair $({\rm q}_A,\ev_A )$ is an adjoint pair in the ordinary sense. Let $M \in ({}_{R_A} \V)_0$, and $X \in [\A, \V]_0$. By the weak and strong  Yoneda Lemma (see \eqref{Yoneda}), we have
$$\V \textrm{-Nat}(\A(A,-)\otimes M, X) \cong \Hom_{ \V_0}(M, [\A, \V](\A(A,-), X))  \cong \Hom_{ \V_0}(M, X(A)). $$
On the other hand,  since the functor $-\otimes_{R_A}-$ is  $\V$-natural in both components and preserves colimits, 
${ \rm q}_A(M)(-)$ is in fact the coequalizer of the following two $\V$-natural transformations
\begin{equation}\label{equalizer:q_A}
\xymatrix{\A(A,-)\otimes R_A \otimes M \ar@<1ex>[rr]^{\circ_{\A}\otimes \id}  \ar@<-1ex>[rr]_{\id \otimes m} && \A(A,-)\otimes M \ar@{=>}[r] & \A(A,-)\otimes_{R_A} M, } 
\end{equation}
which induces the following diagram 
\begin{equation}\label{eq:V-nat}
\xymatrix{
\V\textrm{-Nat}(\A(A,-)\otimes_{R_A} M, X) \ar@{^{(}->}[rr] && \V\textrm{-Nat}(\A(A,-)\otimes M, X) \ar[d]^{\cong}\\
\Hom_{R_A}(M,X(A)) \ar@{^{(}->}[rr] && \Hom_{ \V_0}(M, X(A)) }
\end{equation}
whose top arrow is an injection. 
\sloppy  Given a $\V$-natural transformation $\alpha: \A(A,-)\otimes_{R_A} M \Rightarrow X$, its image in $\Hom_{\V_0}(M, X(A))$ is simply the composition  
$$ \xymatrix{M \ar[r]^-{\cong} & R_A \otimes_{R_A} M \ar[r]^-{\alpha_A} & X(A)}$$ of   $R_A$-morphisms; for the first isomorphism, see \cite[bottom of p. 183]{Par78}. Hence,  there exists an injective mapping  $$\V\textrm{-Nat}(\A(A,-)\otimes_{R_A} M, X) \hookrightarrow \Hom_{R_A} (M,X(A)),$$
which makes the diagram \eqref{eq:V-nat}  commutative.

Conversely,  for a given $R_A$-morphism $f:M \rightarrow X(A)$, the corresponding $\V$-natural transformation $\alpha^f: \A(A,-)\otimes M \Rightarrow X$ is given by the family 
$$\alpha^f_B: \xymatrix{\A(A,B) \otimes M \ar[r]^-{\id \otimes f} & \A(A,B) \otimes X(A) \ar[r]^-{\overline{X}_{AB}} & X(B)}, \quad B \in \mathcal{A},$$
which satisfies
$$\alpha_{B'}^f \circ (\circ_{\A} \otimes \id)=\overline{X}_{BB'} \circ (\id \otimes \alpha_{B}^f): \A(B,B') \otimes \A(A,B) \otimes M \longrightarrow X(B'), \quad B,B' \in \A;$$
see Proposition~\ref{V-natural}.  Since $f$ is an $R_A$-morphism, we have $\overline{X}_{AA} \circ (\id \otimes f)=f \circ m$, so if we take $B=A$, for every $B' \in \A$, the morphism 
$$\alpha_{B'}^f \circ (\circ_{\A} \otimes \id):\ \A(A,B') \otimes R_A\otimes M \longrightarrow X(B') $$
 satisfies
\begin{align*}
\alpha_{B'}^f \circ (\circ_{\A} \otimes \id)& =\overline{X}_{AB'} \circ (\id \otimes \alpha_{A}^f)\\
&=\overline{X}_{AB'} \circ (\id \otimes ( \overline{X}_{AA} \circ (\id \otimes f)  ))\\
&=\overline{X}_{AB'} \circ (\id \otimes  (f \circ m))\\
& = \alpha_{B'}^f \circ (\id \otimes m).
\end{align*}

By the universal property of coequalizers, the $\V$-natural transformation $\alpha^f$ factors through the canonical $\V$-natural transformation given in \eqref{equalizer:q_A}
$$ \A(A,-)\otimes M \Longrightarrow \A(A,-)\otimes_{R_A} M.$$

Thus,  the injective mapping  in \eqref{eq:V-nat} is  in fact surjective. Therefore,
$$ \V\textrm{-Nat}( {\rm q}_A(M), X)  \cong \Hom_{R_A}(M,X(A)).$$ 
It proves that the pair $({\rm q}_A,{\rm ev}_A)$ is an adjoint pair between ordinary categories $({}_{R_A} \V)_0$ and $[\A, \V]_0$.

As for the enriched case, for any object $T$ in $\V$, we have natural isomorphisms
\begin{align*}
 \V\textrm{-Nat}( {\rm q}_A(M),[T, X]) &\cong \Hom_{R_A}(M,[T,X(A)])\\
& \cong \Hom_{R_A}(M \otimes T,X(A))\\
 & \cong \Hom_{ \V_0}(T, {}_{R_A}[M,X(A)]).
\end{align*}
By \eqref{eq:univ_end}, it follows that the end $[\A, \V]({\rm q}_A(M), X )$ exists, and 
 $$[\A, \V]({\rm q}_A(M), X ) \cong {}_{R_A}[M,X(A)].$$
\end{proof}

The standard notions of (co)limits in the enriched category setting are the so-called \textit{weighted} (or \textit{indexed}) (co)limits; see \cite[\S~3.1]{Kelly}. When $[ \A, \V]$ exists, weighted (co)limits in $[\A, \V]$, whenever they exist, are calculated object-wise; see 
\cite[\S~3.3]{Kelly}. Since our focus is on certain ordinary (co)limits (namely, (co)products) in $[\A, \V]_0$, and $[\A, \V]_0$ need not be an ordinary functor category,  we include the following result for completeness.

\subsection{Lemma.}\label{lemma:odinary_limits_functors}
Let $H: \mathcal{C} \longrightarrow [\A, \V]_0$ be a usual functor. The (co)limit of $H$ exists if and only if for every object $A \in \A$, the (co)limit of the induced functor $H_A: \mathcal{C} \rightarrow ( {}_{R_A} \V)_0$, $H_A(C):=H(C)(A)$, exists. As a result, $[\A, \V]_0$ has finite (co)limits.
\begin{proof}
We only prove the statement for the colimit. We already know that, for every object $C \in \mathcal{C}$ and $A \in \A$, $H(C)(A)$ is an $R_A$-left object, so the induced functor $H_A(-):=H(-)(A)$ takes values in $( {}_{R_A} \V)_0$. By Proposition~\ref{prop:adj_q_ev_p},  the functor ${ \rm ev}_A$ preserves (co)limits, which proves  the necessity part of the statement.

 We now prove the sufficiency part. By assumption, for each object $C \in \C$, $H(C): \A \rightarrow \V$ is a $\V$-functor.  So there is   a family of morphisms $\overline{H(C)}_{AB}: \A(A,B)\otimes_{R_A} H(C)(A) \rightarrow H(C)(B)$, $A,B \in \A$, satisfying the composition and unit axioms; see \eqref{morph_X1}. By assumption, for every object $A \in \A$, $\mbox{colim}_{\mathcal{C}} H_A$ exists in $( {}_{R_A} \V)_0$. We let $T(A):= \mbox{colim}_{\mathcal{C}} H_A $. Since $-\otimes_{R_A}-$ is a left adjoint functor, it  preserves colimits, and 
 for every pair of objects $A,B$ in $\A$, we have a natural isomorphism $ \A(A,B) \otimes T(A) \cong \mbox{colim}_{\mathcal{C}} \A(A,B) \otimes_{R_A} H_A$. By the universal property of colimits,  there exists a unique morphism $\mbox{colim}_{\mathcal{C}} \overline{H(-)}_{AB} \colon \A(A,B) \otimes T(A) \cong \mbox{colim}_{\mathcal{C}} \A(A,B) \otimes_{R_A} H_A \longrightarrow T(B)$  satisfying the following commutative diagram
$$
\xymatrix{   
 \A(A,B) \otimes_{R_A} H(C)(A) \ar[rrr]^-{\overline{H(C)}_{AB}} \ar[d] &&&   H(C)(B) \ar[d]\\
 \A(A,B) \otimes_{R_A} T(A) \ar[rrr]_-{\mbox{colim}_{\mathcal{C}} \overline{H(-)}_{AB}} &&& T(B)
}
$$
We let $\overline{T}_{AB}:=\mbox{colim}_{\mathcal{C}} \overline{H(-)}_{AB}$. By the universal property of colimits, the assignment $A\longmapsto T(A)$,  together with the family of morphisms $\overline{T}_{AB}$, $A , B\in \A$, in $\V_0$
satisfies the composition and unit axioms required in Proposition~\ref{X_AB}, and therefore, it gives rise to a $\V$-functor $T$ from $\A$ to $\V$. Moreover, by  the commutative diagram given above and  Proposition~\ref{V-natural}, for each object $C\in \mathcal{C}$,  the canonical morphism $H(C)(A) \longrightarrow T(A)$, $A\in \A$, leads to a $\V$-natural transformation $H(C) \Rightarrow T$. It is a straightforward verification that the $\V$-functor $T: \A \rightarrow \V$ is indeed the colimit of $H$.

As a consequence, since $\V_0$ has finite (co)limits, then so does $[\A, \V]_0$.
\end{proof}

The following result is a  direct consequence of Lemma~\ref{lemma:odinary_limits_functors}.

\subsection{Corollary.}
If $\V_0$ is quasi-abelian, then so is $[\A, \V]_0$.

\subsection{Remark.} From \cite[Theorem~4.2]{HG16}, we already know that if $\V_0$ is abelian, so is $[\A,\V]_0$. Moreover, if $\{G_j\}_{j \in J}$  is a generating set for $\V_0$, then the set 
 $$\{ \A(A,-)\otimes G_j\mid\ j \in J, A \in \A\}$$
 forms a generating set for $[\A, \V]_0$; see \cite[Theorem~4.2]{HG16}. Consequently, if $\V_0$ is Grothendieck, so is $[\A, \V]_0$.
 
\subsection{}\label{q_p_}
For every object $A$ in $\A$, we define
\begin{align*}
{\rm N}^-(A)&:=\{B \in \A\mid\ \A(B,A)\neq 0  \},\\
{\rm N}^+(A)&:=\{B \in \A\mid\ \A(A,B)\neq 0  \}.
\end{align*}

If $\V_0$ admits $  \left| {\rm N}^-(A) \right |$-fold coproducts for every $A\in \A$, then, 
for every $(M_A)_{A \in \A} \in  \prod_{A\in \A } ({}_{R_A} \V)_0$ and every $B \in \A$, the object in $\V_0$ $${ \rm q}((M_A)_{A \in \A})(B):= \bigoplus_{A \in A} { \rm q}_A(M_A)(B)$$
 is well-defined.  By Lemma~\ref{lemma:odinary_limits_functors}, ${ \rm q}((M_A)_{A \in \A})$ is indeed the coproduct of $\V$-functors ${ \rm q}_A(M_A)$, for $A\in \A$, and therefore,  is a $\V$-functor from $\A$ to $\V$. Consequently, this construction defines a well-defined functor 
$${ \rm q}:=\bigoplus_{A \in A} { \rm q}_A \colon \prod_{A\in \A } ({}_{R_A} \V)_0 \longrightarrow [\A, \V]_0.$$

By the universal property of coproducts, it is easy to verify that $( { \rm q}, { \rm ev})$ 
 $$\xymatrix{ \prod_{A\in \A } ({}_{R_A} \V)_0 \ar@/^1em/[rrr]^{  { \rm q}}  &&&    [\A, \V]_0 \ar@/^1em/[lll]^{ \ev} }.
$$
forms an adjoint pair in both the ordinary and the enriched senses. In this case, we let $\xi:{ \rm q} \ev \Longrightarrow \id$ denote the counit of the adjunction.

Similarly, if $\V_0$ admits $ \left| {\rm N}^+(A) \right |$-fold products for every $A\in \A$, we have the following functor   
$${ \rm p}:=\prod_{A \in A} { \rm p}_A \colon \prod_{A\in \A } ({}_{R_A} \V)_0 \longrightarrow [\A, \V]_0,$$
defined by ${ \rm p}((M_A)_{A \in \A}):= \prod_{A \in A} { \rm p}_A(M_A)$, and we have the adjoint  pair $(  {  \ev }, { \rm p} )$
 $$\xymatrix{ \prod_{A\in \A } ({}_{R_A} \V)_0  \ar@/_1em/[rrr]_{{ \rm p}} &&&    [\A, \V]_0 \ar@/_1em/[lll]_{ \ev} }.
$$
In this case, we let  $\eta: \id \Longrightarrow { \rm p} \ev $
denote the unit  of the adjunction.

\subsection{}\label{setup2}
Recall that, up to this point, we have assumed only that $\mathcal{V}_0$ has finite (co)limits. We now additionally assume that $\mathcal{V}_0$ is a weakly idempotent complete additive category. It is immediate that, for any monoid object $R$ in $\V$, $({}_R \V)_0$ is weakly idempotent complete additive, as well. 

\subsection{}\label{lifted_exact}
Suppose that, for every object $A \in \A$, we are given an exact structure  $  \mathcal{E}_A $ on  $  ({}_{R_A}\V)_0$. Since (co)limits in $[\A, \V]_0$ are computed  object-wise (Lemma~\ref{lemma:odinary_limits_functors}), we  can equip $[\A, \V]_0$ with  an exact structure $\mathcal{E}$ for  which each evaluation functor $\ev_A$ at $A$ is an  exact functor. Namely, $\mathcal{E}$    consists of  the   short exact sequences  in $[\A, \V]_0$   such that $\ev_A (\mathbb{E})$ belongs to $\mathcal{E}_A$ for every ${A \in \A}$. We call $\mathcal{E}$ the \textit{induced exact structure} on $[\A, \V]_0$. We recall the following result from \cite[Proposition 4.11]{ECO26}, which is adapted to our setting.

\subsection{Proposition.}\label{prop:proj_inj}
\sloppy Let $ \prod_ {A \in \A} \mathcal{E}_A $ be an  exact structure on the product  category $ \prod_{A\in \A } ({}_{R_A} \V)_0$, and $\mathcal{E}$ be the induced exact structure  on $[\A, \V]_0$. 
For any object $A \in \mathcal{A}$, the following statements hold.
\begin{enumerate}[i)]
\item If $M$ is an  $\mathcal{E}_A$-projective ($\mathcal{E}_A$-injective) object in $({}_{R_A} \V)_0$, then ${ \rm q}_A(M)$ ($ { \rm p}_A(M)$)  is an  $\mathcal{E}$-projective ($\mathcal{E}$-injective) object in $[\A,\V]_0$.

\item Suppose that $\V_0$ admits $ \left| {\rm N}^-(A) \right |$-fold coproducts for every $A\in \A$. If $(\prod_{A\in \A } ({}_{R_A} \V)_0; \prod_ {A \in \A} \mathcal{E}_A)$ has enough projectives,  then so does $( [\A, \V]_0; \mathcal{E})$. Moreover, any $\mathcal{E}$-projective object is a direct summand of ${ \rm q}((M_A)_{A \in \A})$, where $M_A \in \mathcal{E}_A \mbox{-Proj}$ for every $A \in \A$.

\item Suppose that $\V_0$ admits $  \left| {\rm N}^+(A) \right |$-fold products for every $A\in \A$. If $(\prod_{A\in \A } ({}_{R_A} \V)_0; \prod_ {A \in \A} \mathcal{E}_A)$ has enough injectives, then so does $( [\A, \V]_0; \mathcal{E})$. Moreover,   any $\mathcal{E}$-injective object is a direct summand of ${ \rm p}((M_A)_{A \in \A})$, where $M_A \in \mathcal{E}_A \mbox{-Inj}$ for every $A \in \A$.

\end{enumerate}

\subsection{ }\label{define:trivial_exact_str} As a particular case of \eqref{lifted_exact}, for every $A\in \A$, we consider the trivial  exact structure $\mathcal{E}_A^{\text{triv}}$ on $({}_{R_A} \V)_0$. We denote by  $\mathcal{E}_{\textrm{lowt}} $  the induced exact structure on $ [\A, \V]_0$.
In this context, the subscript `lowt' refers to  `locally object-wise trivial'. Applying Proposition~\ref{prop:proj_inj}, we obtain the following.

\subsection{Corollary. }\label{prop:proj_inj_trivial}
\begin{enumerate}[i)]
\item If $\V_0$ admits $   \left| {\rm N}^-(A) \right |$-fold coproducts for every $A\in \A$, then the exact category  $([\A, \V]_0; \mathcal{E}_{\textrm{lowt}})$ has enough projectives. Moreover,  every $\mathcal{E}_{\textrm{lowt}}$-projective  object in $[\A, \V]_0$ is a direct summand of a $\V$-functor of the form
$$ { \rm q} ((M_A)_{A \in \A}), \quad  \quad (M_A)_{A \in \A} \in \prod_{A\in \A} ({}_{R_A} \V)_0.$$

\item If $\V_0$ admits $   \left| {\rm N}^+(A) \right |$-fold products for every $A\in \A$, then the exact category  $([\A, \V]_0; \mathcal{E}_{\textrm{lowt}})$ has  enough injectives.  Moreover, every $\mathcal{E}_{\textrm{lowt}}$-injective object in $[\A, \V]_0$ is a direct summand of a $\V$-functor of the form
$$   { \rm p} ((M_A)_{A \in \A}), \quad  \quad (M_A)_{A \in \A} \in \prod_{A\in \A} ({}_{R_A} \V)_0.$$

\end{enumerate} 


\subsection{Example.}\label{ex:R-enriched-N-chain} Let $R$ be a monoid object in $\V$, and consider the $\V$-category $\A:=\mathcal{L}_{\V} \otimes \mathcal{R}$ given in  Example~\ref{ex_N_chain_N}. As already pointed out, $[\A, \V]_0 \cong  \mathbf{C}_N(({}_R \V)_0)$. For every object $i\in \A$, the locally induced monoid object $R_i$ is $R$, and therefore, the exact structure $\mathcal{E}_{\textrm{lowt}} $ on $[\A, \V]_0$ corresponds to the degreewise trivial exact structure on $\mathbf{C}_N(({}_R \V)_0)$. Since for every object $i $ in $\A$, both  ${\rm N}^-(i)$ and ${\rm N}^+(i)$ are finite sets, $\A$ satisfies both conditions of    Proposition~\ref{prop:proj_inj}. In this case, ${ \rm q}_i$ corresponds to the $i$th disc  functor ${ \rm D}^i$, while  ${ \rm p}_i$ is the $(i+(N-1))$th disc  functor ${ \rm D}^{i+(N-1)}$; see \cite[\S~2.1]{Gil15}.

\subsection{Example.}\label{ex:m_periodic22} Let $R$ be a monoid object in $\V$, and consider the $\V$-category $\A:=\mathcal{L}_{\V} \otimes \mathcal{R}$ given in  Example~\ref{ex_m_periodic11}. As already pointed out, $[\A, \V]_0 \cong  \mathbf{C}_{\equiv_m}(({}_R \V)_0)$. For every object $i\in \A$, the locally induced monoid object $R_i$ is $R$, and therefore, the exact structure $\mathcal{E}_{\textrm{lowt}} $ on $[\A, \V]_0$ corresponds to the degreewise trivial exact structure on $\mathbf{C}_{\equiv_m}(({}_R \V)_0)$. Since  $\A$ has finite objects,  $\A$ satisfies both conditions of    Proposition~\ref{prop:proj_inj}. In this case, ${ \rm q}_i$ corresponds to the $i$th disc  functor ${ \rm D}^i$, while  ${ \rm p}_i$ is the $i+1$th disc  functor ${ \rm D}^{i+1}$.

\subsection{Remark.} Note that one may start with the trivial exact structure  $\mathcal{E}_{\V}^{\text{triv}}$ on $\V_0$. It  induces another `trivial' exact structure, denoted by $\mathcal{E}_{\textrm{owt}}$,  on $[\A, \V]_0$: $\mathcal{E}_{\textrm{owt}}$ consists of   those short exact sequences $\mathbb{E}$ in $[\A, \V]_0$ such that  $\ev_A(\mathbb{E})$ is trivial in $\V_0$ for every $A \in \A$.  The exact structure $\mathcal{E}_{\textrm{owt}}$ can be described in terms of an exact structure on the product category  $ \prod_{A \in \A} ({}_{R_A} \V)_0$. Indeed, for each object $A \in \A$, let  $\mathcal{E}'_A$ denote  the exact structure on $({}_{R_A} \V )_0$  consisting of short exact sequences of $R_A$-left-objects that become trivial after applying the forgetful functor $U: ({}_{R_A} \V )_0 \rightarrow \V_0$. Then
 $$\mathbb{E} \in \mathcal{E}_{\textrm{owt}}\quad  \Longleftrightarrow \quad \ev_A(\mathbb{E}) \in \mathcal{E}_A',  \quad \textrm{ for every } A \in \A.$$
It is immediate that $\mathcal{E}_{\textrm{lowt}} \subseteq  \mathcal{E}_{\textrm{owt}}$. Besides, for every object $A \in \A$,  the  $\mathcal{E}_A'$-projective and $\mathcal{E}_A'$-injective objects in $({}_{R_A} \V )_0$ are precisely the direct summands of $R_A$-left-objects of  the form 
 $$R_A \otimes N \quad \textrm{ and } \quad [R_A, N], $$
 respectively, where $N$ is any object in $\V_0$. Indeed, recall that  $(R_A\otimes-, U, [R_A,-])$ is an adjoint triple between $({}_{R_A} \V)_0$ and $\V_0$ (see \eqref{adj_functors2*} ), and that  the forgetful functor $U$ is exact. It follows that  direct summands of objects of the form $R_A \otimes N$ are $\mathcal{E}_A'$-projective, while direct summands of objects of the form  $[R_A, N]$ are $\mathcal{E}_A'$-injective. Conversely, given $M \in ({}_{R_A} \V)_0$, the counit and unit of the adjunction yield $R_A$-morphisms 
 $$f_m:R_A \otimes U(M) \rightarrow M, \quad \textrm{ and }  \quad g_M:M \rightarrow [R_A, U(M)].$$
 Since $({}_{R_A} \V)_0$ has finite (co)limits, $\Ker(f_M)$ and $\coker (g_M)$ exist. Note that $U(f_M)$ is a retraction, while $U(g_M)$ is a section. Since $\V_0$ is weakly idempotent complete, $U(f_M)$  and $U(g_M)$ are admissible epic and admissible monic, respectively, in $\mathcal{E}_{\V}^{\text{triv}}$. As $U$ creates (co)limits, it follows that $f_M$ is indeed  admissible epic in $\mathcal{E}_A'$, while $g_M$ is  admissible epic in $\mathcal{E}_A'$. Therefore, if $M$ is $\mathcal{E}_A'$-projective, then $f_M$ splits,  so $M$ is a direct summand of $R_A \otimes U(M)$. Similarly,  if $M$ is $\mathcal{E}_A'$-injective, then $g_M$ splits,  so  $M$ is a direct summand of $[R_A, U(M)]$.

\subsection{}\label{morph_phi_psi_X} Let   $X:\A \rightarrow \V$ be a $\V$-functor. If $\V_0$ admits $   \left| {\rm N}^-(A) \right |$-fold coproducts for every $A\in \A$,  then, by the universal properties of coproducts, for every object $A\in \A$,  the $R_A$-morphisms associated with $X$ in \eqref{morph_X1} induce, for each $A\in\mathcal{A}$, a unique $R_A$-morphism
$$\xymatrix{ \varphi_A^X: \displaystyle \bigoplus_{ \substack{ B \in \A \\ A \neq B}} \A(B,A) \otimes_{R_B} X(B) \ar[rr] && X(A)}.$$

Similarly, if $\V_0$ admits $   \left| {\rm N}^+(A) \right |$-fold products for every $A\in \A$,  then, by the universal properties of products, for every object $A\in \A$,  the $R_A$-morphisms associated with $X$ in \eqref{morph_X2} induce, for each $A\in\mathcal{A}$, a unique $R_A$-morphism
$$\xymatrix{  \psi_A^X: X(A) \ar[rr]&& \displaystyle \prod_{ \substack{ B \in \A \\ A \neq B}}{}_{R_B}[\A(A,B),X(B)]}.$$

\subsection{Proposition.}\label{prop:ses1}
Let $X: \A \rightarrow \V$ be a $\V$-functor. 
\begin{enumerate}[i)] 

\item If $\V_0$ admits $  \left| {\rm N}^-(A) \right |$-fold coproducts for every $A\in \A$, then the $\V$-natural transformation
$  \xi_X: { \rm q} \ev (X) \Longrightarrow  X    $
is    and an admissible epic in $\mathcal{E}_{\textrm{lowt}}$.

\item If $\V_0$ admits  $   \left| {\rm N}^+(A) \right |$-fold products for every $A\in \A$, then the $\V$-natural transformation $\eta_X:\ X  \Longrightarrow { \rm p} \ev (X) $ is an admissible monic in $\mathcal{E}_{\textrm{lowt}}$.

\end{enumerate}
\begin{proof}
We only prove the first statement, as similar arguments can be applied to prove the second one.

Note that for every object $A$  in $\A$, the morphism 
$$(\xi_X)_A:{ \rm q} \ev (X)(A)=\bigoplus_{B \in \A}\A(B,A)\otimes_{R_B}X(B) \longrightarrow  X(A)  $$
is in fact
$$
\begin{bmatrix}
\varphi_A^X , \id_{X(A)}
\end{bmatrix}:  \left ( \bigoplus_{\substack{B \in \A \\ A \neq B}}\A(B,A)\otimes_{R_B}X(B) \right ) \oplus X(A) \longrightarrow  X(A)   $$
which is clearly a retraction of  $R_A$-left-objects. So $\xi_X$ is an admissible epic in $\mathcal{E}_{\textrm{lowt}}$. 
 \end{proof}

\subsection{Stable categories of  $[\A, \V]_0$.}\label{subsection_stable}
We denote by  $\mathcal{P}$ and $\mathcal{I}$  the ideals of $\V$-natural transformations in $[\A, \V]_0$ that factor through $\mathcal{E}_{\textrm{lowt}}$-projective and $\mathcal{E}_{\textrm{lowt}}$-injective objects, respectively.  We let 
$$\St_{ \mathcal{P}}([\A, \V]_0) \quad \textrm{ and } \quad \St_{ \mathcal{I}}([\A, \V]_0)$$
 denote the stable  categories of  $[\A, \V]_0$ modulo the ideals $\mathcal{P}$ and $\mathcal{I}$, respectively.

If the condition  in   Proposition~\ref{prop:ses1}-(i) is satisfied, then     every  $\V$-functor $X: \A \rightarrow \V$ admits a short exact sequence in  $\mathcal{E}_{\textrm{lowt}}$ of the form
\begin{equation}\label{eq:canonical.s.e.s1}
\xymatrix{ \Omega(X) \ar@{^{(}->}[r] & { \rm q} \ev (X) \ar@{->>}[r]^-{\xi_X} & X },
\end{equation}
where, by Corollary~\ref{prop:proj_inj_trivial}-(i), ${ \rm q} \ev (X) $ is 
$\mathcal{E}_{\textrm{lowt}}$-projective.  Consequently, by \eqref{Frob},  the stable category   $\St_{ \mathcal{P}}([\A, \V]_0)$ admits a left triangulated structure with the loop functor $\Omega:   \St_{ \mathcal{P}}([\A, \V]_0) \longrightarrow \St_{ \mathcal{P}}([\A, \V]_0)$.

On the other hand, if the condition  in   Proposition~\ref{prop:ses1}-(ii) is satisfied, then     every  $\V$-functor $X: \A \rightarrow \V$ admits a short exact sequence in  $\mathcal{E}_{\textrm{lowt}}$ of the form
\begin{equation}\label{eq:canonical.s.e.s2}
 \xymatrix{  X \ar@{^{(}->}[r]^-{\eta_X} & { \rm p} \ev (X) \ar@{->>}[r] & \Sigma(X)   }.  
 \end{equation}
where, by Corollary~\ref{prop:proj_inj_trivial}-(ii), ${ \rm p} \ev (X) $ is 
$\mathcal{E}_{\textrm{lowt}}$-injective.  Consequently, by \eqref{Frob},  the stable category   $\St_{ \mathcal{I}}([\A, \V]_0)$ admits a right triangulated structure with the suspension functor $\Sigma:  \St_{ \mathcal{I}}([\A, \V]_0)  \longrightarrow \St_{ \mathcal{I}}([\A, \V]_0)$.

In the following two results, we provide explicit descriptions of the $\V$-functors $\Omega(X)$ and $\Sigma(X)$ in the form as given in \eqref{morph_X1} and \eqref{morph_X2}, respectively. 
\subsection{Proposition.}\label{loop_suspension}
Let  $X : \A \rightarrow \V$ be a $\V$-functor. Suppose that $\V_0$ admits  $  \left| {\rm N}^-(A) \right |$-fold coproducts for every $A\in \A$. For any object   $A \in \A$,  we have
$$
\Omega(X)(A)=  \bigoplus_{\substack{B \in \A \\ A \neq B}}\A(B,A)\otimes_{R_B}X(B).
$$
For any pair of  objects   $A, A'\in \A$, the $R_{A'}$-morphism
$$
\overline{\Omega(X)}_{AA'}: \A(A,A') \otimes_{R_A} \Omega(X)(A)    \longrightarrow \Omega(X)(A'),
$$
where $ \A(A,A') \otimes_{R_A} \Omega(X)(A) \cong \bigoplus_{\substack{B \in \A \\ A \neq B}} \A(A,A') \otimes_{R_A}\A(B,A)\otimes_{R_B}X(B)$, is given according to the following two cases:

\begin{itemize}
\item[\underline{Case I:}]   $A=A'$.

$$\overline{\Omega(X)}_{AA}=\bigoplus_{\substack{B \in \A \\ A\neq B}} \circ_{\A}  \otimes \id_{X(B)}.
$$

\item[\underline{Case II:}]  $A\neq A'$.

$$
\overline{\Omega(X)}_{AA'}=(\bigoplus_{\substack{B \in \A \\ A, A'\neq B}} (\circ_{\A} \otimes \id_{X(B)}) )-(\id_{\A(A,A')} \otimes_{R_A} (\varphi_A^X) ), \quad \textrm{ where }
$$

$$
\scalebox{0.8}{\xymatrix{ 
{ \displaystyle \bigoplus_{\substack{B \in \A \\ A \neq B}} \A(A,A') \otimes_{R_A}\A(B,A)\otimes_{R_B}X(B)}  \ar[rrrr]^{-\id_{\A(A,A')} \otimes_{R_A} (\varphi_A^X)} \ar@{->>}[d]_{\pi} &&   && \A(A,A') \otimes_{R_A} X(A)   \ar@{^{(}->}[d]^{\iota}\\
{ \displaystyle \bigoplus_{\substack{B \in \A \\ A, A'\neq B}} \A(A,A') \otimes_{R_A} \A(B,A)\otimes_{R_B} X(B)} \ar[rr]_-{\oplus \circ_{\A} \otimes \id_{X(B)} } && { \displaystyle \bigoplus_{\substack{B \in \A \\ A, A'\neq B}}  \A(B,A') \otimes_{R_B} X(B)} \ar@{^{(}->}[rr]_-{\iota} && \Omega(X)(A')}}
$$
Here, $\pi$ and $\iota$ denote the corresponding canonical projection and inclusion morphisms, respectively.
\end{itemize}

 \begin{proof}
 From the construction of $\xi_X$ in the proof of Proposition~\ref{prop:ses1}, it follows  immediatly that
\begin{align*}
\Omega(X)(A)&=\Ker( \xi_X)(A)=  \bigoplus_{\substack{B \in \A \\ A \neq B}}\A(B,A)\otimes_{R_B}X(B).
\end{align*}
Moreover,   the $\V$-natural transformation $\Omega(X) \Longrightarrow  {\rm q} \ev (X)$ is   given, at any object $A$ in $\A$, by
  \begin{align*}
\begin{bmatrix}
\id, - \varphi_A^X
\end{bmatrix}: \quad &\Omega(X)(A) \longrightarrow {\rm q} \ev (X)(A)= \Omega(X)(A)\oplus X(A). 
\end{align*}
By Proposition \ref{V-natural}, given two objects $A,A'$ in $\A$, we have the following commutative diagrams
$$
\xymatrix{\A(A,A') \otimes_{R_A} \Omega(X)(A) \ar[rrr]^-{ \id \otimes \begin{bmatrix}
\id, - \varphi_A^X
\end{bmatrix}} \ar[d]_{ \overline{\Omega(X)}_{AA'}  } &&& \A(A,A') \otimes_{R_A} (\Omega(X)(A) \oplus X(A)) \ar[d]^{\overline{{ \rm q} \ev (X)}_{AA'}=\oplus (\circ_{\A}\otimes \id)}\\
\Omega(X)(A') \ar[rrr]_-{ \begin{bmatrix}
\id, - \varphi_{A'}^X
\end{bmatrix}} &&&\Omega(X)(A') \oplus X(A')
} 
 $$ 
which yields the description of $\overline{\Omega(X)}_{AA'}$ in each of the respective cases.


 \end{proof}

\subsection{Proposition.}\label{loop_suspension2}
Let  $X : \A \rightarrow \V$ be a $\V$-functor. Suppose that $\V_0$ admits  $  \left| {\rm N}^+(A) \right |$-fold products for every $A\in \A$. For any object   $A \in \A$,  we have
$$\Sigma(X)(A)=  \prod_{\substack{B \in \A \\ A \neq B}} {}_{R_B} [\A(A,B),X(B)].$$

For any pair of  objects   $A, A'\in \A$, the $R_{A}$-morphism
$$
\underline{\Sigma(X)}_{AA'}: \Sigma(X)(A) \longrightarrow {}_{R_{A'}}[\A(A,A'),\Sigma(X)(A')],
$$
where ${}_{R_{A'}}[\A(A,A'),\Sigma(X)(A')] \cong   \prod_{\substack{B \in \A \\ A' \neq B}} {}_{R_B}[ \A(A',B) \otimes_{R_{A'}} \A(A,A'), X(B)   ]$, 
is given according to the following two cases:

\begin{itemize}

\item[\underline{Case I:}]   $A=A'$.

$$
\underline{\Sigma(X)}_{AA}=  \prod_{\substack{B \in \A \\ A\neq B}} {}_{R_B}[ \circ_{\A}, \id_{X(B)}]
$$

\item[\underline{Case II:}]  $A\neq A'$.

$$\underline{\Sigma(X)}_{AA'} =(\prod_{\substack{B \in \A \\ A, A'\neq B}} {}_{R_B}[ \circ_{\A}, \id_{X(B)}])-{}_{R_{A'}}[\id_{\A(A,A')}, \psi_{A'}^X]^*, \quad \textrm{ where }$$

$$
\scalebox{0.7}{\xymatrix{
\Sigma(X)(A) \ar@{->>}[rr]^{\pi} \ar@{->>}[d]_{\pi} && {}_{R_{A'}}[\A(A,A'),X(A')]  \ar[rr]^{-[\id,\psi_{A'}^X]} && {}_{R_{A'}}[\A(A,A'), \Sigma(X)(A')]  \ar[d]^{\cong} \\
{ \displaystyle \prod_{\substack{B \in \A \\ A, A'\neq B}} {}_{R_B}[\A(A,B), X(B)]} \ar[rr]_-{\prod [\circ_{\A}, \id_{X(B)}] } &&  { \displaystyle \prod_{\substack{B \in \A \\ A, A'\neq B}} {}_{R_B}[\A(A',B) \otimes_{R_{A'}} \A(A,A'), X(B)]} \ar@{^{(}->}[rr]_{\iota} && { \displaystyle \prod_{\substack{B \in \A \\ A'\neq B}}  {}_{R_B}[\A(A',B) \otimes_{R_{A'}} \A(A,A'), X(B)]}
}}
$$

Here,  $\pi$ and $\iota$ denote the corresponding canonical projection and inclusion morphisms, respectively; 
${}_{R_{A'}}[\id_{\A(A,A')}, \psi_{A'}^X]^*$ denotes the composition of the morphism 
${}_{R_{A'}}[\id_{\A(A,A')}, \psi_{A'}^X] $
with the canonical isomorphism 
$ {}_{R_{A'}}[\A(A,A'), \Sigma(X)(A')]  \cong { \displaystyle \prod_{\substack{B \in \A \\ A'\neq B}}  {}_{R_B}[\A(A',B) \otimes_{R_{A'}} \A(A,A'), X(B)]} $.

\end{itemize}

 \begin{proof}
 From the constructions of  $\eta_X$ in the proof of Proposition~\ref{prop:ses1}, it follows  immediatly that
\begin{align*}
\Sigma(X)(A)&=\coker(\eta_X)(A)=  \prod_{\substack{B \in \A \\ A \neq B}} {}_{R_B} [\A(A,B),X(B)]. 
\end{align*}
Moreover,   the $\V$-natural transformation ${\rm p} \ev (X) \Longrightarrow \Sigma(X) $ is given, at any object $A$ in $\A$, by
  \begin{align*}
\begin{bmatrix}
\id\\
- \psi_A^X
\end{bmatrix}  : \quad &  {\rm p} \ev (X)(A)= \Sigma(X)(A)\oplus X(A) \longrightarrow \Sigma(X)(A)
\end{align*}
By Proposition \ref{V-natural}, given two objects $A,A'$ in $\A$, we have the following commutative diagram
 $$
\xymatrix{ \Sigma(X)(A) \oplus X(A) \ar[rr]^-{\begin{bmatrix}
\id\\
- \psi_A^X
\end{bmatrix}  } \ar[d]_{\underline{{ \rm p} \ev (X)}_{AA'} }&& \Sigma(X)(A) \ar[d]^{\underline{\Sigma(X)}_{AA'}}\\
[\A(A,A'),\Sigma(X)(A') \oplus X(A')] \ar[rr]_-{ \begin{bmatrix}
\id\\
\begin{bmatrix}
\id\\
-\psi_{A'}^X
\end{bmatrix}
\end{bmatrix} } && [\A(A,A'),\Sigma(X)(A')]} 
 $$
which yield the description of $ \underline{\Sigma(X)}_{AA'}$ in each of the respective cases.


 \end{proof}

 \subsection{Example.}\label{ex:chain1} Let $k:=\mathbb{Z}$ and $R$ be a ring. Consider the category $\A:= \mathcal{L}_{\Z} \otimes  \mathcal{R}$ given in  Example~\ref{ex:R-enriched-N-chain}.
By abuse of notation, for a given $N$-chain   complex $X$ of left $R$-modules, we denote all of its differentials by $d$ 
$$\xymatrix{ \cdots\ar[r]^d & X_{i+1} \ar[r]^d & X_i \ar[r]^d & X_{i-1} \ar[r]^d & \cdots }.$$
By Proposition~\ref{loop_suspension} and \ref{loop_suspension2},   for any integer $n$, we have
\begin{align*}
\Omega(X)_n &=  X_{n+N-1}\oplus \cdots \oplus X_{n+2} \oplus X_{n+1}, \\
\Sigma(X)_n & =X_{n-1}\oplus X_{n-2}\oplus \cdots \oplus X_{n-(N-1)}.
\end{align*}
The   differentials $\Omega(X)_n \rightarrow \Omega(X)_{n-1}$ and $\Sigma(X)_n \rightarrow \Sigma(X)_{n-1}$ are given by 
\begin{align*}
(x_{N-1}, \ldots, x_2, x_1 ) & \mapsto \quad (x_{N-2}, \ldots, x_2, x_1, -d^{N-1}(x_{N-1})-\cdots-d^2 (x_2)-d(x_1) ), \\
(x_1, x_2, \ldots, x_{N-1} ) & \mapsto \quad (x_2-d(x_1),x_3-d^2(x_1), \ldots, x_{N-1}-d^{N-2}(x_1), -d^{N-1}(x_1)), 
\end{align*}
respectively,  where $d^i$ denotes the composition of $i$ consecutive differentials of $X$. We point out that $\Omega(X)$ is precisely the chain complex which  is defined in  the proof of  \cite[Proposition 4.1]{Gil15}. 
   
   In the classical case $N=2$, we recover $\Sigma(X)=X[1]$, which is the suspension functor, while $\Omega(X)=X[-1]$.
   
Analogous element-based descriptions of the loop and suspension functors can be obtained for  $N$-chain complexes in $\mathbf{C}_N(\mathbf{Ban})$.  
\subsection{Example.}\label{ex:m_periodic} Let $m$ be a positive integer, and $R$ be a monoid object in $\V$. Consider the category $\A:=\mathcal{L}\otimes \mathcal{R}$ given in Example~\ref{ex:m_periodic22}. 
For a given $\V$-functor $X:\A \rightarrow \V$, we already know that, for every $i \in \A$,  $X_i:=X(i)$ is an $R$-left-object.   By abuse of notation, for every $i \neq j \in \A$ with $\A(i,j) \neq 0$, we denote by $d$ the  $R$-morphism $X_i \rightarrow X_j$ corresponding, under the tensor-hom adjunction,   to the induced morphism  $R \rightarrow [X_i, X_j]$. So   $X:\A \rightarrow \V$ can be represented by
$$ 
\xymatrix{X_0 \ar[r]^d & X_1 \ar[r]^d \ar[r]& \cdots \ar[r]^d & X_m \ar@/^{5mm}/[lll]^d }.
$$
Applying Proposition~\ref{loop_suspension} and \ref{loop_suspension2}, we have 
\begin{align*}
\Omega(X): \quad  &\xymatrix{X_{m-1} \ar[r]^{-d} & X_0 \ar[r]^{-d} \ar[r]& \cdots \ar[r]^{-d} & X_{m-1} \ar@/^{5mm}/[lll]^{-d} }\\
\Sigma(X):  \quad & \xymatrix{X_1 \ar[r]^{-d} & X_2 \ar[r]^{-d} \ar[r]& \cdots \ar[r]^d & X_0 \ar@/^{5mm}/[lll]^{-d} }.
\end{align*}

\subsection{Example.}\label{ex:morita_context}
Consider a Morita context ring $U:=\begin{bmatrix}
T & M\\
N & S
\end{bmatrix}$, that is,   $T$, $S$ are rings, $M$ and $N$ are $T$-$S$ and $S$-$T$-bimodules, respectively, together with  bimodule homomorphisms
$M \otimes_S N \rightarrow T$ , and $N \otimes_T M \rightarrow S$.

 We let $\A$ be the category with two objects $\{ 0,1 \}$ and  with morphisms $\Hom(0,0):=T$, $\Hom(1,1):=S$,  $\Hom(0,1):=N$, $\Hom(1,0):=M$ together with the canonical composition rule induced from the ring multiplications and bimodule homomorphisms of the given Morita context ring.
 It is well known  that
$$\Add(\A, \Ab) \cong U \Mod. $$
Note that a left   $U$-module $X$ is of the form $(X_0,X_1, \varphi_0, \varphi_1)$, where $X_0$ and $X_1$ are  left $T$- and $S$-modules, respectively, and  $\varphi_0: M \otimes_S X_1 \rightarrow X_0$ and $\varphi_1: N \otimes_T X_0 \rightarrow X_1$ are $T$- and $S$-module homomorphisms, respectively, subject to associativity axioms. The left $U$-module $\Omega(X)$ is of the form $(M \otimes_S X_1,N \otimes_T X_0, -\id_M \otimes \varphi_1, -\id_N \otimes\varphi_0 )$, while $\Sigma(X)$ is
$$(\Hom_S(N,X_1), \Hom_T(M,X_0), -\varphi_1^*,-\varphi_0^*),$$
where  $\varphi_1^*$ and $\varphi_1^*$ are the following  canonical compositions, respectively,
$$ M\otimes_S \Hom_T(M,X_0) \longrightarrow X_0 \longrightarrow \Hom_S(N,X_1),$$
$$ N\otimes_T \Hom_S(N,X_1) \longrightarrow X_1 \longrightarrow \Hom_T(M,X_0).$$
Here, the first homomorphisms are the counit morphisms of the tensor-hom adjunction, while the second homomorphisms correspond, under the tensor-hom adjunction,  to $\varphi_1$ and $\varphi_0$, respectively.

\vspace{5mm}

We conclude this section by adressing the following question.

\subsection{Question.}\label{question} Is there a $\V$-category  whose underlying ordinary category is $\St_{ \mathcal{P}}([\A, \V]_0)$, or $\St_{ \mathcal{I}}([\A, \V]_0)$?

\vspace{5mm}

The answer to this question is affirmative under certain conditions, and we leave its details to the appendix.
 
 \section{Homotopy in functor categories}\label{section:homotopy}
Throughout this section, we continue to assume that $\V_0$ is a weakly idempotent complete additive category with finite (co)limits. In this section, we further study  the stable categories $\St_{ \mathcal{P}}([\A, \V]_0)$ and $\St_{ \mathcal{I}}([\A, \V]_0)$ (see \eqref{subsection_stable}).  We first establish sufficient conditions on $\mathcal{A}$ under which $([\A, \V]_0; \mathcal{E}_{\textrm{lowt}})$ forms a Frobenius exact category; that is, when $\mathcal{E}_{\textrm{lowt}}$-projective and $\mathcal{E}_{\textrm{lowt}}$-injective objects in $[\A, \V]_0$ coincide. Subsequently, we consider the conditions under which  
left or right homotopic to zero $\V$-natural transformations in $[\A, \V]_0$ can be characterized in terms of a family of $R_A$-morphisms, analogous to the case of  chain complexes.  For it, we introduce the following  weak and enriched versions of Serre functors. 

\subsection{Definition.}\label{def:left_right_weak_Serre}
A \textit{right-weak Serre $\V$-functor on $\A$} is a $\V$-functor ${ \rm F}: \A \longrightarrow \A$ together with a  family $\sigma=\{ \sigma_A \}_{A \in \A}$  of $\V$-natural transformations in $[\A, \V]_0$
$$\sigma_A:{}_{R_{{\rm F}(A)}}[ \A(-,  {\rm F}(A)), R_{{\rm F}(A)}] \Longrightarrow \A(A,-), \quad {A \in \A},$$
such that for every $B\in \A$, $(\sigma_A)_B$ is a retraction as $R_B$-$R_A$-biobjects.\footnote{The $\V$-functor $F$ induces a monoid morphism $R_A \longrightarrow R_{F(A)}$, and therefore, the $R_A$-right-object structure on ${}_{R_{{\rm F}(A)}}[ \A(-,  {\rm F}(A)), R_{{\rm F}(A)}]$ is induced  from the $R_{F(A)}$-$R_A$-biobject structure on $R_{F(A)}$.}\footnote{As already pointed out in \eqref{evaluation}, $(\sigma_A)_B$ is already an $R_B$-morphism. So being $R_B$-$R_A$-biobject retraction means that $ (\sigma_A)_B$ is a retraction as an $R_A$-morphism with a right inverse which is $R_B$-morphism.}

A \textit{left-weak Serre $\V$-functor on $\A$} is  a  $\V$-functor ${ \rm G}: \A \longrightarrow \A$ together with a family $\nu=\{ \nu_A \}_{A \in \A}$ of $\V$-natural transformations in $[\A^{\textrm{op}}, \V]_0$
$$\nu_A: [ \A( {\rm G}(A),-), R_{{\rm G}(A)}]_{R_{{\rm G}(A)}} \Longrightarrow \A(-,A), \quad {A \in \A},$$
such that for every $B \in \A$, $(\nu_A)_B$ is a retraction as $R_A$-$R_B$-biobjects.

\subsection{Remark.}\label{remark:duality_left_serre} It is straightforward to observe that if  $({\rm F }, \sigma)$ is a right-weak Serre  $\V$-functor on $\A$, then $({\rm F}^{\textrm{op}}, \sigma)$ is a left-weak Serre $\V$-functor on $\A^{\textrm{op}}$.

\subsection{Remark.}\label{remark:serre_R_extension} One advantage of right/left-weak $\V$-functors is that they are stable under the $R$-extension of $\A$ for any monoid object $R$. More precisely, if $(F, \sigma)$ is a right-weak Serre $\V$-functor on $\A$, then it  induces naturally a right-weak Serre $\V$-functor on the $R$-extension $\V$-category $\A \otimes \mathcal{R}$. Indeed, for every objects $A \in \A$, we have 
$${}_{R_{{\rm F}(A)} \otimes R}[ \A(-,  {\rm F}(A)) \otimes R, R_{{\rm F}(A)} \otimes R]={}_{R_{{\rm F}(A)}}[ \A(-,  {\rm F}(A)), R_{{\rm F}(A)}] \otimes R,$$
and the functor $- \otimes R$ preserves retractions. Dually, a left-weak Serre $\V$-functor $(G, \nu)$ on $\A$,   induces naturally a left-weak Serre $\V$-functor on the $R$-extension $\V$-category $\A \otimes \mathcal{R}$.

\subsection{Remark (A comparison with the classical Serre functor).}\label{Remark:classical_Serre} Let $k$ be a commutative ring, and $\V:=k \Mod$. Let 
$\A$ be a $\Hom$-finite $k$-linear category, that is, $\Hom(A,B)$ is a finitely generated projective $k$-module for every objects $A$ and $B$ in $\A$. The classical Serre functor on  $\A$ relative to $k$  (see \cite[Definition 4.1]{DSS17}) is an autoequivalence $F: \A \rightarrow \A$ together with a natural isomorphism between $k$-linear bifunctors
$$\sigma: \Hom_k(\Hom(-,F(?)),k) \longrightarrow \Hom(?,-).$$

Although  the notions introduced in Definition~\ref{def:left_right_weak_Serre}  have the adjective `weak', there is, in general, no a priori implication from the existence of classical Serre functors and that of weak right or left Serre functors. Namely, a $k$-lineal category may have a classical Serre functor  ${\rm F}$ relative to $k$ without $F$ being   either a right-weak  or a left weak Serre functor  in the sense of Definition~\ref{def:left_right_weak_Serre}, since the latter notions involve locally induced endomorphism rings. However, we have the implication for two cases. 

The first case is when  $R_A=k$ for every object $A\in \A$.  In this situation, the formula given in the classical Serre functor relative to $k$ coincides with those in Definition~\ref{def:left_right_weak_Serre}, and therefore, a classical Serre functor $(F, \sigma)$ induces naturally both a right-weak and a left-weak Serre  $k$-linear functor on $\A$. By Remark~\ref{remark:serre_R_extension}, for any $k$-algebra $R$, they also induce  a right-weak and a left-weak Serre  $k$-linear functor on its $R$-extension category $\A \otimes \mathcal{R}$. A  prototypical example of this case is the free $k$-linear category $\mathcal{L}_k$ considered in Example~\ref{ex:enriched-Qco}, which   is  known to admit a classical Serre functor defined by $F(i):=i-(N-1)$; see \cite[Paragraph~2.4]{HJ19}. Consequently, $\mathcal{L}_k \otimes R$ admits a right-weak and left-weak Serre  $k$-linear functor. In Example~\ref{ex:N_complex_Serre}, we extend this construction to a more general setting, allowing an arbitrary $\V$ and a monoid object $R$ in $\V$. In the same \cite[Paragraph~2.4]{HJ19}, it is 
 noted   that  the  $k$-linear category with $R=k$ considered in Example~\ref{ex:m_periodic} also admits a classical Serre functor defined by $F(i):=i+1$ for $0 \leq i <m-1$, and $F(m-1):=0$. So for any $k$-algebra $R$, its $R$-extension category naturally admits both left/right-weak Serre $k$-linear functor; see Example~\ref{Serre_m_periodic}. We point out that, in both cases, the corresponding $R$-extension $k$-linear categories cannot admit a classical Serre functor relative to $k$ unless $R$ is a symmetric $k$-algebra.

The second case in which the existence of a classical Serre functor relative to $k$ implies the existence of right-weak and left-weak Serre functors is when $F=\id$. Indeed, suppose that the $k$-linear category $\A$ has a classical Serre functor $(\id, \sigma)$ relative to $k$. By the naturality of $\sigma$ in both components, for every object $A \in \A$, the $k$-linear isomorphism induced from $\sigma$
$$R_A^*:=\Hom_k(R_A,k) \longrightarrow R_A,$$ 
is an $R_A$-$R_A$-bimodule isomorphism, and therefore, $R_A$ is a symmetric $k$-algebra; see \cite[p. 49]{CB25}. So for any $B \in \A$, we have
\begin{align*}
\Hom_{R_A}(\Hom(B,A),R_A) & \cong \Hom_{R_A}(\Hom(B,A),R_A^*)\\
& \cong \Hom_k(R_A \otimes_{R_A} \Hom(B,A),k)\\
& \cong \Hom_k( \Hom(B,A),k)\\
& \cong \Hom(A,B),
\end{align*}
where the last isomorphism is induced from $\sigma$. As an example, take a finite dimensional symmetric $k$-algebra $R$ (for instance, $R:=k(x)/{(x^2)}$) , and consider the category $\A$ with two objects corresponding to the Morita context ring $\begin{bmatrix}
R & R \\
R & R
\end{bmatrix} $; see Example~\ref{ex:morita_context}.
\subsection{Example.}\label{ex:N_complex_Serre} Let $R$ be a monoid object in $\V$. Consider the  $\V$-category $ \A:=\mathcal{L}_{\V} \otimes \mathcal{R}$  given in Example~\ref{ex_N_chain_N}.
 We define  a $\V$-functor ${\rm F}: \A \longrightarrow \A $ as follows: For each pair of objects $i,j \in \A$, set ${\rm F}(i):=i-(N-1)$,   and let the induced morphism on hom-objects  $\A(i,j) \rightarrow \A(F(i), {\rm F}(j))$ be  the identity morphism. It is straightforward to verify that ${\rm F}$ is indeed a $\V$-functor. In fact, ${\rm F}$ is a $\V$-autoequivalence with the inverse ${\rm G}(i)=i+(N-1)$. Note that for every object $i \in \A $, $R_{i}=R$, and we define 
\begin{align*}
\sigma_i&:=\{ (\sigma_i)_j \}_{j \in \Z}: {}_R[ \A(-, F(i)) , R] \Longrightarrow \A(i,-)\\
\nu_i & :=\{ (\nu_i)_j \}_{j \in \Z}: [\A(G(i),-), R]_R \Longrightarrow \A(-, i)
\end{align*}
as follows: For any integer $j$, 
\begin{align*}
\A(i,j)=R \Longleftrightarrow \A(i,j) \neq 0 &\Longleftrightarrow i-(N-1) \leq j \leq i\\
& \Longleftrightarrow j-(N-1) \leq i-(N-1) \leq j\\
& \Longleftrightarrow \A (j, {\rm F}(i)) \neq 0\\
& \Longleftrightarrow \A(j, {\rm F}(i))=R
\end{align*}
Whenever these equivalent conditions hold, we  let $(\sigma_i)_j$ be the canonical isomorphism ${}_R[R,R] \longrightarrow R$. Otherwise, we take $(\sigma_i)_j$ to be  the zero morphism.

Similarly, we define $(\nu_i)_j$ by considering the following equivalences
\begin{align*}
\A(j,i)=R \Longleftrightarrow \A(j,i) \neq 0 &\Longleftrightarrow i \leq j \leq i+(N-1)\\
& \Longleftrightarrow \A ({\rm G}(i),j) \neq 0\\
& \Longleftrightarrow \A({\rm G}(j),i)=I.
\end{align*}

\subsection{Example.}\label{Serre_m_periodic}
Let $R$ be a monoid object in $\V$. Consider the $\V$-category given in  Example~\ref{ex:m_periodic}. By similar arguments applied in Example~\ref{ex:N_complex_Serre}, it can be easily verified that the canonical $\V$-functors on $\A$ defined by 
$$F(i):=i+1, \quad  \textrm{ for }  \quad 0 \leq i < m-1,  \quad \textrm{ and  } F(m-1):=0, $$
$$G(i):=i-1, \quad  \textrm{ for }  \quad 0 < i  \leq  m-1,  \quad \textrm{ and  } F(0):=m-1$$
are parts of natural right-weak and left-weak Serre $\V$-functors on $\A$, respectively.

\vspace{5mm}

 For simplicity,  we provide a set of examples in their $k$-linear enriched form.  However, the same results hold for their $\V$-enriched counterparts for an arbitrary $\V$.

\subsection{Example.}
\label{ex:non_Serre} Let $k$ be a commutative ring, and $R$ be a non-Frobenius $k$-algebra\footnote{$R$ is not isomorphic to $\Hom_k(R,k)$ as a left $R$-module.}(see (\cite[p.~49]{CB25}) (for ex., the upper triangular $2\times 2$ matrix ring over $k$). Let $S:=k \times R$ be the product ring. We consider the Morita context ring
$U:=\begin{bmatrix}
k & M\\
N & S
\end{bmatrix}$,
where $M:=k$ and $N:=k$. The right and left $S$-scalar multiplications on $M$ and $N$,respectively,  are induced by  the canonical projection ring homomorphism $\pi: S \rightarrow k$. The $k$-$S$ and $S$-$k$-bimodule homomorphisms, respectively, are given by
\begin{align*}
\varphi:& M \otimes_S N \rightarrow k, \quad &\varphi(m  \otimes n):=m\cdot n,\\
\psi: & N \otimes_k M \rightarrow S, \quad &\psi(n \otimes m):=(nm, 0_R).
\end{align*}
We let $\A$ denote the $k$-linear category with two objects $\{ 0, 1\}$ corresponding to the Morita context ring $U$ as described in  Example~\ref{ex:morita_context}. We claim that the identity functor $\id: \A \longrightarrow \A$ is a part of both a left-weak and a right-weak k-linear Serre functor on $\A$. First, we define the right-weak $k$-linear Serre structure, that is,  two natural transformations $\sigma_i$, $i=0,1$, where 
$$\sigma_i:\Hom_{R_i} ( \Hom(-,  i), R_i) \Longrightarrow \Hom(i,-), \quad i=0,1,$$
such that for every $j\in \A$, $(\sigma_i)_j$ is a retraction of  $R_j$-$R_i$-bimodules. We let $\sigma_0$ denote the family consisting of the following $k$-linear homomorphisms
$$(\sigma_0)_0=(\sigma_0)_1:\Hom_k(k,k) \longrightarrow k, \quad f \rightsquigarrow f(1),$$
where $\Hom(1,0)=M=k=N=\Hom(0,1)$.
We first show that $\sigma_0$ is a natural transformation. It is  equivalent to showing that $(\sigma_0)_1$ is a left $S$-module homomorphism, and that, for every $n \in N$ and $m \in M$, the following diagrams commute
$$
\scalebox{1}{\xymatrix{ f \ar@{~>}[d] \\
f( \varphi(- \otimes n))
}}
\scalebox{1}{\xymatrix{ \Hom_k(k,k)\ar[r]^-{(\sigma_0)_0} \ar[d] & k \ar[d]^{n\cdot(-)}\\
 \Hom_k(M,k) \ar[r]_-{(\sigma_0)_1} & N}}
 \quad 
 \scalebox{1}{\xymatrix{f \ar@{~>}[d] \\
 f(-\cdot m) }}
\scalebox{1}{\xymatrix{ \Hom_k(M,k) \ar[r]^-{(\sigma_0)_1} \ar[d] & N \ar[d]^{\varphi(m \otimes -)}\\
 \Hom_k(k,k) \ar[r]_-{(\sigma_0)_0} & k}}
$$
The commutativity of the previous diagrams is immediate. It remains to show that $(\sigma_0)_1$ is a left $S$-module homomorphism. For it,  let $(\alpha, r) \in S $, and $f\in \Hom_k(M,k)$. Then we have
$$(\sigma_0)_1( (\alpha, r)  \cdot f  )=(\sigma_0)_1 ( f(-\cdot (\alpha, r)  ))= (\sigma_0)_1 ( f(-\cdot \alpha  ))=f(\alpha)=\alpha f(1)=(\alpha, r)\cdot f(1).$$
Finally,  it is immediate that  $(\sigma_0)_0$ is an isomorphism of  $k$-modules, and $(\sigma_0)_1$ is an isomorphism  of $S$-$k$-bimodules. 
 

Now, we let $\sigma_1$ denote the family consisting of the following $k$-linear homomorphisms
$$  (\sigma_1)_0=\Hom_S(N,S) \longrightarrow M, \quad f \rightsquigarrow \pi_1(f(1)), $$
$$(\sigma_1)_1=\Hom_S(S,S) \longrightarrow S, \quad f \rightsquigarrow f(1),$$
where $\pi_1: S \rightarrow k$ is the canonical projection mapping. Both homomorphisms are clearly $k$-module homomorphisms, and $(\sigma_1)_1$ is an $S$-$S$-bimodule isomorphism. So  the naturality of $\sigma_1$ is equivalent to showing that for every $n \in N$ and $m \in M$, the following diagrams commute 
$$
\xymatrix{f \ar@{~>}[d] \\
f( - \cdot  n)}
\scalebox{1}{\xymatrix{ \Hom_S(N,S)\ar[r]^-{(\sigma_1)_0} \ar[d] & M \ar[d]^{\psi(n\otimes(-))}\\
 \Hom_S(S,S) \ar[r]_-{(\sigma_1)_1} & S}}
 \quad
 \xymatrix{f \ar@{~>}[d]\\
 f(\psi(- \otimes m))}
\scalebox{1}{\xymatrix{\Hom_S(S,S) \ar[r]^-{(\sigma_1)_1} \ar[d] & S \ar[d]^{m \cdot(-)}\\
 \Hom_S(N,S) \ar[r]_-{(\sigma_0)_0} & M}}
$$
Note that if $f:N \rightarrow S$ is an $S$-module homomorphism with $(\alpha, r):=f(1)$, then
$$(\alpha, r)=f(1)=f((1,0_R) \cdot 1)=(1, 0_R) \cdot f(1)=(\alpha,0).$$
So for every $n \in N$, $f(n)=f((n, 0_R)\cdot 1)=(n \alpha,0_R)$
By definition, $(\sigma_1 )_0(f)=\alpha$, and therefore, the commutativity of the first diagram follows. For the commutativity of the second diagram, let $f:S \rightarrow S$ be an $S$-module homomorphism with $f(1,1)=(\alpha, r)$. Then for any $m \in M$, we have 
$$ \pi_1(f (\psi (1 \otimes m)))=\pi_1( f(m, 0_R))=\pi_1( (m,0)\cdot (\alpha,r))=m \cdot \alpha=m \cdot (\alpha, r).$$ 
It remains to show that $(\sigma_1)_0$ is a $k$-$S$-bimodule retraction. It is clearly a right $S$-module homomorphism. Recall that the right $S$-module structure on $\Hom_S(N, S)$ is induced from $S_S$. As pointed out before, a left $S$-module homomorphism $f: N \rightarrow S$ is completely determined by an unique scalar $\alpha \in k$. So $(\sigma_1)_0$ is a bijection. 

Similar constructions and arguments can be applied to show that $\id:\A^{\textrm{op}}  \rightarrow \A^{\textrm{op}} $ is a part of a right-weak $k$-linear Serre functor, where  $\A^{\textrm{op}}$ corresponds to the Morita context ring $\begin{bmatrix}
S & N\\
M & k
\end{bmatrix}$. By Remark~\ref{remark:duality_left_serre} ,  this implies the existence of left-weak k-linear Serre functor on $\A$.

\subsection{Remark.}\label{remark:explan_non_Serre} There are two important points to note regarding Example~\ref{ex:non_Serre}. First, there do not exist a $k$-algebra $R' \neq k$ and a $k$-linear category $\mathcal{A}'$ such that $\mathcal{A}' \otimes \mathcal{R}' = \mathcal{A}$ (see \eqref{R-ext}).  
Indeed, if such a $k$-algebra $R'$ and a $k$-linear category $\mathcal{A}'$ existed, then we would have a $k$-algebra isomorphism $\End_{\mathcal{A}'}(0)\otimes_k R' \cong k$, which would imply $\End_{\mathcal{A}'} (0)\cong k \cong R'$ as an algebra. Second, even if we assume that the non-Frobenius $k$-algebra $R$ is  finite dimensional, $\mathcal{A}$ doesn't a classical Serre functor $(F, \sigma)$ relative to $k$. Indeed,  if such a Serre functor existed, then, by the naturality of  $\sigma$, the isomorphism $\Hom_k(\Hom(1,F(1)),k)\cong \Hom(1,1)=S$ of $k$-modules would be an isomorphism of left $S$-modules. So if $F(1)=0$, then we would have  $k \cong \Hom_k(\Hom(1,0),k) \cong S=\Hom(1,1)$ as $S$-$k$-bimodules. However,  $k$ is not isomorphic to $S$ as a left $S$-module. So we would   have $F(1)=1$. In this case,  again by the naturality of $\sigma$, we have an isomorphism of left $S$-modules $\Hom_k(S,k)=\Hom_k(\Hom(1,1),k)\cong  \Hom(1,1)=S$. Since $S=k \times R$, we have an equivalence of module categories $S \Mod \cong (k\Mod) \times (R\Mod)$, and therefore, $\Hom_k(\Hom(1,1),k)\cong \Hom_k(k,k) \times \Hom_k(R,k) \cong S= k \times R$, which implies $\Hom_k(R,k) \cong R$ as left $R$-modules, which is a contradiction by the assumption of being non-Frobenius for $R$. As noted at the beginning, the failure to admit a classical Serre functor relative to \(k\) in this case is due neither to being $R'$-extension of  another category nor to the requirement that $F$ be an autoequivalence. Rather, it stems from the requirement that $\sigma$ is a natural isomorphism in both variables, which is restrictive.

\subsection{Example.}\label{ex:right_repetitive_quiver} Let $n \geq 2$, and $R$ be a ring. We consider a slightly modified version of the repetitive quiver $\mathbb{Z} { \rm A}_n$ with mesh relations presented in \cite[Example 8.10]{HJ22}. We define the quiver $\mathbb{Z^+} A_n$ with mesh relations as follows: ${ \rm A}_n \times \mathbb{Z}^+$ is the set of vertices, and  the set of arrows are formed by 
$$a_{i,j}: (i,j) \rightarrow (i+1, j), \quad \textrm{ and } \quad a_{i,j}^*: (i+1,j) \rightarrow (i, j+1),$$
where  $1 \leq  i \leq n-1$ and $j \in \mathbb{Z}^+$, together with the mesh relations
$$
0=a^*_{1,j} a_{1,j}=a_{n-1,j+1} a_{n-1,j}^*=a_{i+1,j}^*a_{i+1,j}+a_{i,j+1} a_{i,j}, \quad j\geq 1, 1\leq i \leq n-2.
$$
For $n=3$, the quiver ${ \rm A}_3 \times \mathbb{Z}^+$ is of the form
$$
\scalebox{0.9}{\xymatrix{
(1,1) \ar[rd]^{a_{1,1}} &  & (1,2) \ar[rd]^{a_{1,2}} &  & (1,3)\ar[rd]^{a_{1,3}} & & (1,4) \ar[rd]^{a_{1,4}} &  \\
 & (2,1) \ar[ru]^{a^*_{1,1}}  \ar[rd]^{a_{2,1}} & & (2,2) \ar[rd]^{a_{2,2}}  \ar[ru]^{a^*_{1,2}}  & & (2,3) \ar[ru]^{a^*_{1,3}}  \ar[rd]^{a_{2,3}} & & \ldots \\
& & (3,1) \ar[ru]^{a^*_{2,1}}  && (3,2) \ar[ru]^{a^*_{2,2}} && (3,3)\ar[ru]^{a^*_{2,3}} & 
}}
$$

We denote by  $\A$ the  $R$-extension of the path category   of the quiver $\mathbb{Z^+} A_3$  modulo mesh relations. Note that, when $n=2$, we have $  \Add(\A, \Ab) \cong \bfC_-(R) $, where $\bfC_-(R)$ is  the category of non-positive chain complexes of left $R$-modules.

The category $\A$ does not admit a Serre functor relative to $\mathbb{Z}$. In fact, suppose that there exists
a Serre  functor $({ \rm F}, \sigma)$.  The  natural transformation $\sigma: \Hom_{\Z}(\Hom(-{ \rm F}(?)),\Z) \cong \Hom(?,-)$ would force ${ \rm F}$ to have
$${ \rm F}(1,j)=(n,j), { \rm F}(i,j)=(i,j+1), { \rm F}(n,j)=(1,j+2), \quad j \geq 1, 1 < i < n.$$
However, ${ \rm F}$ can not be an autoequivalence, since ${ \rm F}$ is not essentially surjective. 
Indeed, there is no nonzero morphism
in $\A$ from any of $(i,1)$, $1\leq i<n$, and $(1,2)$ to an object in
$\operatorname{im}({ \rm F})$. On the other hand, ${ \rm F}$ naturally forms part of a
right-weak Serre functor on $\A$. In fact, 
$$\Hom((i,j),(k,l))\neq 0 \quad \Leftrightarrow \quad \Hom((k,l),{ \rm F}(i,j)) \neq 0.$$
Note that $\Hom((i,j),(k,l))\neq 0$ if and only if $\Hom((i,j),(k,l))=R$.

\subsection{Example.}\label{ex:right-chain} Let $R$ be a ring, and   $\A$ be the category whose objects are positive integers,
and  whose morphisms are given by 
$$
\Hom(m,n):= \begin{cases}   
R,& \textrm{ if } m \leq n \leq 2m\\
0, & \textrm{ otherwise} 
\end{cases}
$$
together with the ring  multiplication as the composition rule in $\A$.
We define an additive functor ${ \rm F}: \A \rightarrow \A$ by ${ \rm F}(m):=2m$ together with the identity morphism on morphisms.
Let $m \geq 1$ be  a fixed integer. We define a natural transformation $$\sigma_m := \{ (\sigma_m)_n \}_{n 
\geq 1}: \Hom_R(\Hom(-,2m), R) \Longrightarrow \Hom(m,-)$$ 
as follows: For a given integer $n \geq 1$, note that $\Hom(m,n)\neq 0$ if and only if $m \leq n \leq 2m$ if and only if $n \leq 2m \leq 2n$ if and only if $\Hom(n, { \rm F} (m)) \neq 0$. Under these conditions, we define $(\sigma_m)_n$ to be  the canonical isomorphism 
 $$ \Hom_R(\Hom(n,2m),R)= \Hom_R(R,R) \cong R =\Hom(m, n).$$ 
Otherwise, we set 
 $(\sigma_m)_n$ to be the zero morphism 
$$  \Hom_R(\Hom(n,2m),R)=0=\Hom(m, n).$$ 
One can easily verify that the  family $\sigma_m := \{ (\sigma_m)_n \}_{n 
\geq 1}$ is in fact a natural isomorphism.
So 
${ \rm F}$ together with $\{ \sigma_m \}_{m \geq 1}$ is a right-weak Serre additive functor on $\A$.


\vspace{5mm}

Now, we present  two lists of  conditions that will be needed to prove Theorem~\ref{theo:proj-inj}.

\subsection{Condition.}\label{conditions:frob}
\begin{enumerate}
\item 
For every object $A$ in $\A$, the set ${\rm N}^-(A)$ is finite.
\item 
 For every pair of objects $A$ and $B$ in $\A$, $\A(A,B)$ is a dualizable   $R_B$-left-object (see \eqref{internally_proj}).
\item $\A$ has a right-weak Serre functor.
\end{enumerate}

\subsection{Condition.}\label{conditions:frob2}
\begin{enumerate}
\item 
For every object $A$ in $\A$, the set
 ${\rm N}^+(A)$ is finite
\item 
 For every pair of objects $A$ and $B$ in $\A$, $\A(A,B)$ is a  dualizable $R_A$-right-object. 
 
\item $\A$ has a left-weak Serre functor.
\end{enumerate}

\subsection{Theorem.}\label{theo:proj-inj}
\begin{enumerate}[i)]
\item If Condition~\ref{conditions:frob} is satisfied, then  $\mathcal{E}_{\textrm{lowt}} \mbox{-Proj} \subseteq  \mathcal{E}_{\textrm{lowt}} \mbox{-Inj}$.
\item If Condition~\ref{conditions:frob2} is satisfied, then  $ \mathcal{E}_{\textrm{lowt}} \mbox{-Inj} \subseteq  \mathcal{E}_{\textrm{lowt}} \mbox{-Proj}$.

\end{enumerate}
\begin{proof}
\begin{enumerate}[i)]
\item 
From Corollary~\ref{prop:proj_inj_trivial}, we already know that any $\mathcal{E}_{\textrm{lowt}}$-projective object is a direct summand of a $\V$-functor of the form 
$$ { \rm q} ((M_A)_{A \in \A}), \quad (M_A)_{A \in \A} \in \prod_{A \in \A} ({}_{R_A} \V)_0.$$  
By assumption, for every object $B \in \A$,   ${\rm N}^-(B)$ is finite, and therefore,  
$$ { \rm q} ((M_A)_{A \in \A})(B)=\bigoplus_{A \in \A} \A (A, B) \otimes_{R_A} M_A  = \prod_{A \in \A} \A (A, B) \otimes_{R_A} M_A,$$
that is, 
$${ \rm q} ((M_A)_{A \in \A})=\bigoplus_{A \in \A}{\rm q}_A(M_A) \cong \prod_{A \in \A}{\rm q}_A(M_A). $$
Thus,  it is enough to show that for every object $A$ in $\A$ and every $R_A$-left-object $M$, ${\rm q}_A(M)$  is    $\mathcal{E}_{\textrm{lowt}}$-injective.  By assumption, there exists a right-weak Serre  $\V$-functor $({ \rm F}, \sigma)$ on $\A$. Recall that  for every object $B$ in $\A$ $(\sigma_A)_B$ is a retraction as $R_B$-$R_A$-biobjects, so we have the following $\V$-natural transformation of $\V$-functors from $\A$ to $\V$
$$\sigma_A \otimes \id_M: {}_{R_{{\rm F}(A)}}[\A(-,{ \rm F}(A)),{R_{{\rm F}(A)}}]\otimes_{R_A} M \Longrightarrow \A(A,-)\otimes_{R_A} M= {\rm q}_A(M),  $$
which is  an admissible epic in  $\mathcal{E}_{\textrm{lowt}}$.  Here, $R_{{ \rm F}(A)}$ is considered to be an $R_{{ \rm F}(A)}$-$R_A$-biobject. Additionally, as $\A(-,{ \rm F}(A))$ is a dualizable   ${R_{{\rm F}(A)}}$-left-object, by Proposition~\ref{prop:dualiz_S_R}, 
\begin{align*}
{}_{R_{{\rm F}(A)}}[\A(-,{ \rm F}(A)),{R_{{\rm F}(A)}}]\otimes_{R_A} M &\cong {}_{R_{{\rm F}(A)}}[\A(-,{ \rm F}(A)),{R_{{\rm F}(A)}} \otimes_{R_A} M]\\
&= {\rm p}_{{\rm F}(A)}( {R_{{\rm F}(A)}} \otimes_{R_A} M), 
\end{align*}
and therefore,  we have an admissible epic in $\mathcal{E}_{\textrm{lowt}}$ 
$${\rm p}_{{\rm F}(A)}( {R_{{\rm F}(A)}} \otimes_{R_A} M) \Longrightarrow {\rm q}_A(M)$$
As ${\rm q}_A(M)$ is an $\mathcal{E}_{\textrm{lowt}}$-projective object, it splits in $[\A, \V]_0$, and hence, ${\rm q}_A(M)$ is a direct summand of ${\rm p}_{{\rm F}(A)}( {R_{{\rm F}(A)}} \otimes_{R_A} M) $, which is an $\mathcal{E}_{\textrm{lowt}}$-injective object by  Corollary~\ref{prop:proj_inj_trivial}.

\item 
From Corollary~\ref{prop:proj_inj_trivial}, we know that any $\mathcal{E}_{\textrm{lowt}}$-injective object is a direct summand of a $\V$-functor of the form 
 $$ { \rm p} ((M_A)_{A \in \A}), \quad (M_A)_{A \in \A} \in \prod_{A \in \A} ({}_{R_A} \V)_0.$$ 
 As for every object $B \in \A$, ${\rm N}^+(B)$ is finite, 
$$ { \rm p} ((M_A)_{A \in \A})(B)=\prod_{A \in \A} {}_{R_A} [\A (B,A),  M_A]  = \bigoplus_{A \in \A} {}_{R_A} [\A(B,A),  M_A].$$ 
So we have
$${ \rm p} ((M_A)_{A \in \A})=\prod_{A \in \A}{\rm p}_A(M_A) \cong \bigoplus_{A \in \A}{\rm p}_A(M_A). $$
Thus, it is enough to show that ${ \rm p}_A(M) $ is $\mathcal{E}_{\textrm{lowt}}$-projective for any $M \in ({}_{R_A}\V)_0$.  By assumption, there exists a left-weak Serre functor $( { \rm G}, \nu)$ on $\A$. As for every object $B$ in $\A$, $(\nu_A)_B$ is a retraction as $R_A$-$R_B$-biobjects, we have 
the following $\V$-natural transformation ${}_{R_A}[\nu_A,M]$ in $[\A, \V]_0$ 
$${ \rm p}_A(M)= {}_{R_A}[\A(-,A),M] \Longrightarrow  {}_{R_A}[ [ \A( {\rm G}(A),-), R_{{\rm G}(A)}]_{R_{{\rm G}(A)}},M], $$ 
which is in fact an admissible monic in $ \mathcal{E}_{\textrm{lowt}}$, and therefore, splits. By assumption,  $\A( { \rm G}(A),-)$ is a dualizable $R_{ { \rm G}(A)}$-right-object, by Proposition~\ref{prop:dualizable_codualizable},  there is a $\V$-natural isomorphism 
\begin{align*}
 {}_{R_A} [ [\A({\rm G}(A),-), R_{{\rm G}(A)}]_{R_{{\rm G}(A)}},M] \cong  &\A({ \rm G}(A),-) \otimes_{R_{G(A)}} {}_{R_A}[R_{G(A)},M]\\
 &={\rm q_{{{ \rm G}(A)} }}({}_{R_A}[R_{G(A)},M])
 \end{align*}
So  ${\rm p}_A(M)$ is a direct summand of ${\rm q_{ {{ \rm G}(A)} }}( {}_{R_A}[R_{G(A)},M]) $, which is $\mathcal{E}_{\textrm{lowt}}$-projective by Corollary~\ref{prop:proj_inj_trivial}.

\end{enumerate}
\end{proof}

\subsection{Example.} The categories given in Examples~\ref{ex:right_repetitive_quiver} and \ref{ex:right-chain}  satisfy Conditions~\ref{conditions:frob}, while their opposite categories verify Conditions~\ref{conditions:frob2} .

\subsection{Corollary.}\label{thm:frob} If Conditions~\ref{conditions:frob} and \ref{conditions:frob2} are satisfied, then $([\A, \V]_0; \mathcal{E}_{\textrm{lowt}})$ is a Frobenius exact category.

\subsection{Example.} Recall that any monoid object $R$ in $\V$ is a dualizable object both in ${}_R \V$ and $\V_R$. Hence the $\V$-categories given in Examples~\ref{ex:N_complex_Serre} and \ref{Serre_m_periodic} satisfy both Conditions~\ref{conditions:frob} and \ref{conditions:frob2}.

\subsection{Example.}\label{ex:non_Serre2}
Consider $k$-linear category $\A$ given in Example~\ref{ex:non_Serre}.  Clearly, $k$ is a finitely generated projective as left and right $S$-module, since it is direct summand of $S$. By Example~\ref{ex:dual_module}, it is a dualizable left and right $S$-module. So both Conditions~\ref{conditions:frob} and \ref{conditions:frob2} are satisfied. By Corollary~\ref{thm:frob}, $(k \operatorname{-Lin}(\A, \Ab); \mathcal{E}_{\textrm{lowt}})$ is a Frobenius exact category.

\vspace{5mm}

Under the conditions given in Proposition~\ref{prop:ses1} (i) and (ii), a $\V$-natural transformation $\alpha: X \Longrightarrow Y$ in $[\A, \V]_0$ is left homotopic to zero if and only if  it factors through $\xi_Y: { \rm q ev}(Y) \Longrightarrow Y$, and is right homotopic to zero if and only if it factors through $\eta_X: X \Longrightarrow { \rm p ev}(X)$, respectively. This way of formulating  left/right homotopic to zero is somewhat abstract. In the following result, we provide an equivalent characterization in terms of a family of morphisms in $\V_0$.

\subsection{Theorem.}\label{thm:homotopy} 
\begin{enumerate}[i)]
\item Suppose that $\A$ satisfies Condition~\ref{conditions:frob2}. A $\V$-natural transformation $\alpha:X \Longrightarrow Y$ in $[\A, \V]_0$ is right homotopic to zero if and only if there exists a family $\{s_A\}_{A \in \A}$ of 
 $R_{{ \rm G}(A)}$-morphisms 
$$s_{A}: {}_{R_A}[R_{{ \rm G}(A)}, X(A)]  \longrightarrow Y({ \rm G}(A))$$  such that $\alpha$ is equal to the composition 
$$\xymatrix{X \ar@{=>}[r]^-{  \eta_X} & { \rm p} { \rm ev}(X) \ar@{=>}[rr]^-{({}_{R_A}[\nu_A, X(A)])}&&  \bigoplus_{A \in \A} { \rm q}_{{ \rm G}(A)}({}_{R_A}[R_{{ \rm G}(A)}, X(A)]) \ar@{=>}[rr]^-{\Sigma_{A \in \A} \overline{\beta}_A} && Y},$$
where $\overline{\beta}_A$ denotes  the $\V$-natural transformation corresponding to $s_A$ under the adjunction $({ \rm q}_{{ \rm G}(A)}, { \rm ev}_A)$
$$\Hom_{R_{{ \rm G}(A)}} ( {}_{R_A}[R_{{ \rm G}(A)}, X(A)], Y ({ \rm G}(A))) \cong \V \mbox{-Nat}({ \rm q}_{{ \rm G}(A)}({}_{R_A}[R_{{ \rm G}(A)}, X(A))], Y ).$$
\item Suppose that $\A$ satisfies Condition~\ref{conditions:frob}. A $\V$-natural transformation $\alpha:X \Longrightarrow Y$ in $[\A, \V]_0$ is left homotopic to zero if and only if there exists a family $\{s_A\}_{A \in \A}$ of 
 $R_{{ \rm F}(A)}$-morphisms 
$$s_{A}:X({ \rm F}(A)) \longrightarrow  R_{{ \rm F}(A)} \otimes_{R_A} Y(A)  $$  such that $\alpha$ is equal to the composition 
$$\xymatrix{X \ar@{=>}[rr]^-{ ( \overline{\beta}_A)_{A \in \A}} && \prod_{A \in \A} { \rm p}_{F(A)}(R_{{ \rm F}(A)} \otimes_{R_A} Y(A)) \ar@{=>}[rr]^-{(\sigma_A \otimes Y(A))_{A\in \A}}&&  { \rm q} { \rm ev}(Y)  \ar@{=>}[r]^-{\xi_Y} & Y},$$
where $\overline{\beta}_A$ denotes  the $\V$-natural transformation corresponding to $s_A$ under the adjunction $({ \rm ev}_{{ \rm F}(A)}, { \rm p}_{{ \rm F}(A)})$
$$\Hom_{R_{{ \rm F}(A)}} ( X({ \rm F}(A)), R_{{ \rm F}(A)} \otimes_{R_A} Y(A)) \cong \V \mbox{-Nat}(X, { \rm p}_{{ \rm F}(A)}(R_{F(A)} \otimes_{R_A} Y(A)) ).$$
\end{enumerate}
\begin{proof}
We prove only the first statement, as the second follows by dual arguments. Since ,$\operatorname{N}^+(A)$ is finite for every object $A$ in $\A$, we have
$${\rm p} \ev (X)=\prod_{A \in \A} {\rm p}_A(X(A))= \bigoplus _{A \in \A} {\rm p}_A(X(A)). $$
So a $\V$-natural transformation $\alpha:X \Longrightarrow Y $ is right homotopic to zero if and only if, for every object $A \in \A$, there exists a $\V$-natural transformation $\beta_A: {\rm p}_A(X(A)) \Longrightarrow Y$ such that $\alpha=(\Sigma_{A \in \A} \beta_A) \circ \eta_X$.

From the proof of Theorem~\ref{theo:proj-inj}, we know that  ${\rm p}_A(X(A))$ is a direct summand of ${ \rm q}_{ { \rm G}(A)}({}_{R_A}[R_{ { \rm G}(A)}, X(A)])$ via the $\V$-natural transformation ${}_{R_A}[\nu_A, X(A)]$. 
So for every object $A \in \A$, there exists a $\V$-natural transformation $\beta_A: {\rm p}_A(X(A)) \Longrightarrow Y$ such that $\alpha=(\Sigma_{A \in \A} \beta_A) \circ \eta_X$ if and only if, for every object $A \in \A$, there exists a $\V$-natural transformation $\overline{\beta}_A: { \rm q}_{{ \rm G}(A)}({}_{R_A}[R_{ { \rm G} (A)}, X(A)])  \Longrightarrow Y$ such that 
$$\alpha=(\Sigma_{A \in \A} \overline{\beta}_A \circ {}_{R_A}[\nu_A, X(A)]) \circ \eta_X.$$

The statement follows from the adjunction given in Proposition~\ref{prop:adj_q_ev_p}
\begin{align*}
\V\mbox{-Nat}({ \rm q}_{ { \rm G} (A)}({}_{R_A}[R_{ { \rm G} (A)}, X(A)]), Y ) & \cong  \Hom_{R_{{ \rm G}(A)}} ( {}_{R_A}[R_{{ \rm G}(A)}, X(A)], Y ({ \rm G}(A))).
\end{align*}
\end{proof}

\subsection{}\label{particular_case_homotopy}Even though  Theorem~\ref{thm:homotopy} explains the role of families of morphisms in describing homotopies between chain complexes, its formulation still remains somewhat abstract in this general setting. To make the result even more concrete, we consider a particular case. Suppose that for every $A \in \A$, we have  $R_A=R$, where $R$ is a fixed monoid object in $\V$, and  that for every two objects $A \neq B$ in $\A$, $\A(A,B)=\bigoplus_{t \in T_{AB}} R_t$ for some finite set $T_{AB}$, where $R_t=R$ for every $t$. For every element $t\in T_{AB}$, $e_t^{AB}: R_t \hookrightarrow \bigoplus_{i \in T_{AB}} R_t$  denotes the canonical inclusion in $\V_0$, which is clearly an $R$-$R$-bimorphism. Note that the indexed set $T_{AB}$ may be  empty which would imply  $\A(A,B)=0$. Under these assumptions, we now consider the following two  cases separately.

\subsubsection{}
Suppose, in addition,  that $  \A$ satisfies   Condition~\ref{conditions:frob2}. Without loss of generality, we continue to use the notation $e_t^{BA}$ for the canonical inclusion of $R_t \cong {}_R[ R_t,R] $ in ${}_R[\A(B,A),R]$, where 
$${}_R[\A(B,A),R]= {}_R[\bigoplus_{t \in T_{BA}} R_t,R] \cong \bigoplus_{t \in T_{BA}} {}_R[ R_t,R]   \cong \bigoplus_{t \in T_{BA}} R_t=\A(B,A).$$  
By assumption, for every object $A \in \A$, there exists a $\V$-natural transformation
$$\nu_A: [\A( { \rm G}(A),-),R]_R \cong \bigoplus_{t \in T_{ { \rm G}(A)\star }} R_t \Longrightarrow \A(-,A)= \bigoplus_{t \in T_{  \star A}} R_t$$
such that for every $B\in \A$, $(\nu_A)_B$ is a retraction as $R$-$R$-biobjects. Suppose that, for every object $B \in \A$, the matrix representation of $(\nu_A)_B$ has entries consisting of zeros and identity morphisms, with a unique identity morphism in each column. Equivalently, it is induced by a surjective mapping
$$
\nu_{BA}\colon T_{ { \rm G}(A)B} \longrightarrow T_{BA}.
$$
such that, for every $t \in T_{{ \rm G}(A)B}$, the following diagram commutes
$$
\xymatrix{R_t \ar@{=}[rr] \ar@{^{(}->}[d]_{e_t^{ { \rm G}(A)B}} && R_{\nu_{BA} (t)}\ar@{^{(}->}[d]^{e_{\nu_{BA} (t) }^{BA}}\\
{ \displaystyle \bigoplus_{t \in T_{ { \rm G} (A)B}} R_t} \ar[rr]_{(\nu_A)_B} && { \displaystyle \bigoplus_{t \in T_{BA}} R_t}
}
$$

Under these conditions,
\begin{align*}
{\rm p}_A(X(A))& \cong \bigoplus_{t\in T_{\star A}}X(A)_t,
\\
{\rm q}_{{ \rm G}(A)} ( {}_{R_{{ \rm G}(A)}} [R_{{ \rm G}(A)},X(A)])
& \cong {\rm q}_{{ \rm G}(A)}(X(A))
\cong \bigoplus_{t\in T_{{ \rm G}(A)\star}}X(A)_t,
\qquad X(A)_t=X(A).
\end{align*}

Furthermore, the section
$
{}_{R_A}[\nu_A,X(A)]\colon
{\rm p}_A(X(A))\Longrightarrow {\rm q}_{ { \rm G}(A)}(X(A))
$
is, at each object $B$, given by the commutative diagram
$$
\xymatrix{X(A)_t \ar[rr]^{\Delta_{X(A)_t}} \ar@{^{(}->}[d]_{e_t^{BA}} && { \displaystyle \bigoplus_{ \substack{ j \in T_{ { \rm G} (A)B}\\ \nu_{BA}(j)=t}} X(A)_j }\ar@{^{(}->}[d]^-{e_{\nu_{BA} (t) }^{BA}}\\
{ \displaystyle \bigoplus_{t \in T_{BA}} X(A)_t} \ar[rr]_{(\nu_A)_B} && { \displaystyle \bigoplus_{j \in T_{ { \rm G} (A)B}} X(A)_j}
}
$$
where $\Delta_{X(A)_t}$ denotes the diagonal morphism of $X(A)_t$, and the rightmost vertical morphism is the canonical inclusion.

Under these conditions, a $\V$-natural transformation $\alpha :X \Longrightarrow Y$  in $[\A, \V]_0$ is right homotopic to zero if and only if for every $A \in \A$, there exists an $R$-morphism $s_A:X(A) \longrightarrow Y(G(A))$ such that  for every $B \in \A$, the $R$-morphism $\alpha_B:X(B) \longrightarrow Y(B)$ is equal to the composition
$$
\xymatrix{
X(B) \ar[r]^-{(\eta_X)_B} &  { \displaystyle \bigoplus_{A \in \A} (\bigoplus_{t \in T_{BA}} X(A)_t)}    \ar[rrr]^{\oplus_{A \in \A} ( \oplus_{t \in T_{BA}} \Delta_{X(A)_t)}}  &&&  { \displaystyle \bigoplus_{A \in \A} (\bigoplus_{t \in T_{ { \rm G} (A)B}} X(A)_t)}  \\
\ar[r]^-{\oplus_{t \in {T_{ { \rm G} (A)B}}}s_A}&{ \displaystyle \bigoplus_{A \in \A} (\bigoplus_{t \in T_{{ \rm G}(A)B}} Y( { \rm G}(A))_t)}    \ar[rr]^-{\oplus_{A\in \A}(\overline{Y}_{ { \rm G} (A)B})} && Y(B) 
}
$$


\subsubsection{}
Suppose, in addition, that $  \A$ satisfies  Condition~\ref{conditions:frob}. 
By assumption, for every object $A \in \A$, there exists a $\V$-natural transformation
$$\sigma_A: {}_R [\A(-, F(A)),R] \cong \bigoplus_{t \in T_{\star F(A)}} R_t \Longrightarrow \A(A,-)= \bigoplus_{t \in T_{  A \star }} R_t$$
such that for every $B\in \A$, $(\sigma_A)_B$ is a retraction as $R$-$R$-biobjects. Suppose that, for every object $B \in \A$, the matrix representation of $(\sigma_A)_B$ has entries consisting of zeros and identity morphisms, with a unique identity morphism in each column. Equivalently, it is induced by a surjective mapping
$$
\sigma_{AB}\colon T_{BF(A)} \longrightarrow T_{AB}.
$$
such that, for every $t \in T_{BF(A)}$, the following diagram commutes
$$
\xymatrix{R_t \ar@{=}[rr] \ar@{^{(}->}[d]_{e_t^{BF(A)}} && R_{\sigma_{AB} (t)}\ar@{^{(}->}[d]^{e_{\sigma_{AB} (t) }^{AB}}\\
{ \displaystyle \bigoplus_{t \in T_{BF(A)}} R_t} \ar[rr]_{(\sigma_A)_B} && { \displaystyle \bigoplus_{t \in T_{AB}} R_t}
}
$$

Under these conditions,
\begin{align*}
{\rm q}_A(Y(A))& \cong \bigoplus_{t\in T_{A\star }}Y(A)_t,
\\
{\rm p}_{F(A)} ( R_{F(A)} \otimes_{R_A} Y(A))
& \cong {\rm p}_{F(A)}(Y(A))
\cong \bigoplus_{t\in T_{\star F(A)}}Y(A)_t,
\qquad Y(A)_t=Y(A).
\end{align*}
Furthermore, the retraction
$
\sigma_A \otimes_R Y(A)\colon
{\rm p}_{F(A)}(Y(A))\Longrightarrow {\rm q}_{A}(X(A))
$
is, at each object $B$, given by the commutative diagram
$$
\xymatrix{Y(A)_t \ar@{=}[rr] \ar@{^{(}->}[d]_{e_t^{BF(A)}} && Y(A)_{\sigma_{AB} (t)}\ar@{^{(}->}[d]^{e_{\sigma_{AB} (t) }^{AB}}\\
{ \displaystyle \bigoplus_{t \in T_{BF(A)}} Y(A)_t} \ar[rr]_{(\sigma_A)_B} && { \displaystyle \bigoplus_{t \in T_{AB}} Y(A)_t}
}
$$

Under these conditions, a $\V$-natural transformation $\alpha :X \Longrightarrow Y$  in $[\A, \V]_0$ is left homotopic to zero if and only if for every $A \in \A$, there exists an $R$-morphism $s_A:X(F(A)) \longrightarrow Y(A)$  such that  for every $B \in \A$, the $R$-morphism $\alpha_B:X(B) \longrightarrow Y(B)$ is equal to the composition
$$
\scalebox{0.9}{\xymatrix{
X(B) \ar[rr]^-{\oplus_{A \in \A} \underline{X}_{BF(A)}} && { \displaystyle \bigoplus_{A \in \A} ( \bigoplus_{t \in T_{B F(A)}} X(F(A))_t) } \ar[rrr]^-{\oplus_{A \in \A} \oplus_{t \in T_{BF(A)}} s_A} &&& { \displaystyle \bigoplus_{A \in \A} ( \bigoplus_{t \in T_{B F(A)}} Y(A)_t) }  \\
 \ar[rr]^-{ \oplus_{A \in \A} \sigma_A \otimes_R Y(A)} &&  { \displaystyle \bigoplus_{A \in \A} ( \bigoplus_{t \in T_{A B}} Y(A)_t) } \ar[rr]^{(\xi_A)_B} &&Y(B)& }}
$$

\subsection{Example.} Let $R$ be a monoid object in $\V$. Consider Example \ref{ex:N_complex_Serre}. A morphism $\alpha:X \rightarrow Y$ in $[  \mathcal{L}_{\V}  \otimes \mathcal{R}, \V]_0 \cong \mathbf{C}_N(R)$ is  homotopic to zero if and only if for every integer $n$, there exists a $R$-morphism $s_n: X_n \longrightarrow X_{n+(N-1)}$ such that 
$$\alpha_n= \sum_{i=0}^{N-1}  d_Y^{(N-1)-1} \circ s_{n-i} \circ d_X^i.$$
So $\St([  \mathcal{L}_{\V}  \otimes \mathcal{R}, \V]_0) \cong \mathbf{K}_N(R)$, the usual $N$-homotopy category of $R$-left-objects.

\subsection{Example.}\label{ex:chain_varying}
Let $R$ be a monoid object in $\V$. Consider Example~\ref{ex:right-chain}. Without loss of generality, for a given $\V$-functor $X$  in $[\A, \V]_0$,  we let $d_X$ denote the $R$-morphism $X_n \rightarrow X_{n-1}$ for any $n \in \Z$, and  $d^m_X$ denote the composite of $m$ consecutive $R$-morphisms $d_X$. 
A $\V$-natural transformation  $\alpha:X \longrightarrow Y$ in $[\A, \V]_0$ is left homotopic to zero if and only if there exists a family $\{ s_n:X_{2n} \longrightarrow Y_{n}\}_{n \geq 1}$ of $R$-morphisms such that 

$$
\alpha_n=\begin{cases}
 \sum_{k=0}^{\frac{n}{2}} d_Y^{\frac{n}{2}-k  } \circ s_{\frac{n}{2} +k} \circ d_X^{2k} & n \textrm{ is even,}  \\
 & \\
 \sum_{k=0}^{\frac{n-1}{2}} d_Y^{\frac{n-1}{2} + k  } \circ s_{\frac{n+1}{2} +k} \circ d_X^{2k+1}  & n \textrm{ is odd}
\end{cases}
$$

$$
\xymatrix{   X_{1} \ar[r]^{d_X} \ar[d]_{\alpha_{1}} & X_{2} \ar[r]^{d_X} \ar@{-->}[dl]_{s_{1}} \ar[d]_{\alpha_{2}}  & X_3 \ar[r]^{d_X}  \ar[d]_{\alpha_{3}}& X_{4} \ar[r]^{d_X}  \ar@{-->}[dll]^<<<<<<<{s_{2}} \ar[d]_{\alpha_{4}}& X_{5}  \ar[r]^{d_X} \ar[d]_{\alpha_{5}} & X_{6} \ar[r]^{d_X}   \ar[d]_{\alpha_{6}} \ar@{-->}[dlll]_{s_{3}} & X_{7} \ar[r] \ar[d]_{\alpha_{7}} & \\
 Y_{1} \ar[r]_{d_Y} & Y_{2} \ar[r]_{d_Y} & Y_3 \ar[r]_{d_Y} & Y_{4} \ar[r]_{d_Y} & Y_{5} \ar[r]_{d_Y} & Y_{6} \ar[r]_{d_Y} & Y_{7} \ar[r] & }
$$

\section{Appendix: Internal $\V$-ideals}

In this section, we address the Question~\ref{question} by  constructing  a $\V$-category whose underlying ordinary category is $\St_{\mathcal{P}}( [ \A, \V]_0)$. 
Throughout, we assume that $\V_0$ satisfies the following setting.

\subsection{Setup.}\label{setup3} $\V_0$ is a weakly idempotent complete quasi-abelian category with finite (co)limit, and  admits  $   \left| {\rm N}^-(A) \right |$-fold coproducts for every $A\in \A$. We further assume that $[\A,\V]$ exists, and that the unit object $I$ is $\mathcal{E}_{\textrm{max}}$-projective in $\V_0$ (Proposition~\ref{prop:quasi_abelian}).

\vspace{3mm}

For a  given $\V$-functor $Y: \A \rightarrow \V$
we consider the short exact sequence given in \eqref{eq:canonical.s.e.s1} for $Y$.
Since ${ \rm q} \ev (Y)$ is $\mathcal{E}_{\textrm{lowt}}$-projective  (Corollary~\ref{prop:proj_inj_trivial}-(i)), a $\V$-natural transformation from $X$ to $Y$ belongs to $\mathcal{P}(X,Y)$ if and only if it factors through $\xi_Y $. 
Applying the $\V$-functor $[\A, \V](X,-)$, which is left exact,  we obtain the following sequence in $\V_0$
$$\xymatrix{ [\A, \V](X, \Omega(Y)) \ar@{^{(}->}[r] &  [\A,\V](X, { \rm q} \ev(Y)) \ar[rr]^-{[\A, \V](X, \xi_Y)} && [\A,\V](X, Y) }.$$
Since $\V_0$ is quasi-abelian, the morphism $[\A, \V](X, \xi)$ 
has the following canonical factorization 
 \begin{equation}\label{eq:factor_q.a.}
 \xymatrix{[\A,\V](X, { \rm q} \ev(Y)) \ar[r] & \im ([\A, \V](X, \xi_Y)) \ar@{^{(}->}[r]^-{\iota} & [\A,\V](X, Y)}  
 \end{equation}
where $\im ([\A, \V](X, \xi_Y)):=\Ker(\coker ( [\A, \V](X, \xi_Y)))$; see \eqref{eq:analysis_qa}.  Since the composition of epics is epic, by Proposition~\ref{prop:epi_monic_q.a}, the first morphism is epic, and 
\begin{equation}\label{eq_E_max_ses}
\xymatrix{ \im ([\A, \V](X, \xi_Y)) \ar@{^{(}->}[r]^-{\iota} & [ \A, \V](X,Y) \ar@{->>}[r] & \coker(\iota)=\coker ( [\A, \V](X, \xi_Y))}.
\end{equation}
is a short exact sequence in $(\V_0; \mathcal{E}_{\textrm{max}})$.  Now, we define
$$
\mathcal{P}_{\V}(X,Y):= \im ([\A, \V](X, \xi_Y)).
$$


 \subsection{Proposition.}\label{prop:enriched_ideal}
 $\mathcal{P}_{\V}(-,-): [\A, \V]^{\textrm{op}} \times [\A, \V] \longrightarrow \V$ is a $\V$-subfunctor of $[ \A, \V](-,-).$

\begin{proof}
We only prove that for any $\V$-functor $X: \A \rightarrow \V$, $\mathcal{P}_{\V}(X,-): [\A, \V] \longrightarrow \V$ is a $\V$-subfunctor of $[ \A, \V](X,-)$.

Consider the $\V$-natural transformation $[\A,\V](X, \xi_{(-)}): [\A,\V](X, { \rm q} \ev(-)) \Rightarrow [\A,\V](X,-) $.  By  Proposition~\ref{V-natural}, for any pair of $\V$-functors $Y,Y' \in [\A, \V]$,
we have the following commutative diagram
$$
\scalebox{0.73}{\xymatrix{
[\A,\V](Y,Y') \otimes [\A,\V](X, { \rm q} \ev(Y)) \ar[rr]^-{\id \otimes [\id,\xi_Y]} \ar[d] && [\A,\V](Y,Y') \otimes [\A,\V](X,Y) \ar@{->>}[rr] \ar[d]^-{\circ} && [\A,\V](Y,Y') \otimes \coker (\iota) \ar@{.>}[d]\\
[\A,\V](X, { \rm q} \ev(Y')) \ar[rr]_{[\id,\xi_{Y'}]} && [\A,\V](X, Y') \ar@{->>}[rr] && \coker([\id, \xi_{Y'}])
}
}
$$
Note that the upper right horizontal morphism is the cokernel of $\id \otimes \iota$, since $\otimes$ preserves cokernels. By the universal property of cokernels, the rightmost vertical morphism exists and is unique. Moreover, since $\otimes$ preserves epics, the canonical morphism $[\A,\V](Y,Y') \otimes [\A,\V](X, { \rm q} \ev(Y)) \rightarrow [\A,\V](Y,Y') \otimes \mathcal{P}_{\V}(X,Y)$ induced by \eqref{eq:factor_q.a.} is still epic. By the universal property of kernels, there exists a unique morphism $\overline{\mathcal{P}_{\V}^X }_{YY'}: [\A,\V](Y,Y') \otimes \mathcal{P}_{\V}(X,Y) \longrightarrow  \mathcal{P}_{\V}(X,Y')$ which makes the following diagram commutative
 $$
\xymatrix{
[\A,\V](Y,Y') \otimes \mathcal{P}_{\V}(X,Y)  \ar[rr]^-{\overline{\mathcal{P}_{\V}^X }_{YY'}} \ar[d]_{\id \otimes \iota } && \mathcal{P}_{\V}(X,Y')\ar@{^{(}->}[d]^{\iota}\\
[\A,\V](Y,Y') \otimes [\A,\V](X,Y)  \ar[rr]^-{\circ } && [\A, \V](X, Y').
} 
 $$
 By the associativity of the composition and the unit axiom in $[\A, \V]$, it is straightforward to verify that the assignment $\mathcal{P}_{\V}(X,Y)$ together with the family of morphisms $\overline{\mathcal{P}_{\V}^X }_{YY'}$, for  $Y, Y' \in [\A, \V]$, satisfies the composition and the unit axioms required in Proposition~\ref{X_AB}, and therefore, $\mathcal{P}_{\V}(X,-): \A \rightarrow \V$ is a $\V$-functor. Furthermore, by Proposition~\ref{V-natural}, the commutative diagram above shows that the inclusion $ \mathcal{P}_{\V}(X,-) \hookrightarrow [\A, \V](X,-)$ is a $\V$-natural transformation.

\end{proof}

 The following result follows immediately by using the fact that $\otimes$ preserves cokernels and  epics.
 \subsection{Proposition.}\label{prop:enriched_quotient} The collection of all $\V$-functors from $\A$ to $\V$, together with the $\V$-hom object
 \begin{align*}
 \underline{[\A, \V]}(X,Y)&:= \coker ( [\A, \V](X, \xi_Y))=\coker(\iota)
  \end{align*}
where $\iota: \mathcal{P}_{\V}(X,Y) \rightarrow [\A, \V](X,Y)$, along with the canonical identity morphism and the composition rule, forms a  $\V$-category, denoted by $\underline{[\A, \V]}$.  Besides, its underlying ordinary category is  precisely the stable category of $[\A, \V]_0$ modulo $\mathcal{P}$, that is, 
 $$(\underline{[\A, \V]})_0=\St_{ \mathcal{P}}([\A, \V]_0) $$
 
 \begin{proof}
 Using the fact that $\otimes $ preserves cokernels and epics, the composition rule in $[\A, \V]$ induces the composition rule in  $\underline{[\A, \V]}(X,Y)$. It is  straightforward to verify that these structures satisfy the axioms of a $\V$-category.

On the other hand, recall that $\Hom_{\V_0}(I, [\A, \V](X, Y))=\V \mbox{-Nat}(X, Y)$.  Since $I$ is an $\mathcal{E}_{\textrm{max}}$-projective object in $\V_0$, applying $\Hom_{\V_0}(I,-)$ to the  short exact sequence in \eqref{eq_E_max_ses}, 
we have the following short exact sequence of abelian groups
$$\scalebox{0.85}{\xymatrix{0 \ar[r] & \Hom_{\V_0}(I, \mathcal{P}_{\V}(X,Y)) \ar[r] & \Hom_{\V_0}(I, [\A,\V](X,Y)) \ar[r] & \Hom_{\V_0}(I,  \underline{[\A, \V]}(X,Y))  \ar[r] & 0  }}   $$
and $\Hom_{\V_0}(I, \mathcal{P}_{\V}(X,Y))$ is the image of $\Hom_{\V_0}(I, [\A, \V](X, \xi_Y))=\V \mbox{-Nat}(X, \xi_Y)$, which is precisely $\mathcal{P}(X,Y)$. Moreover, 
\begin{align*}
\V \mbox{-Nat} (X,Y)/ \mathcal{P}(X,Y) &=\Hom_{\V_0}(I, [\A, \V](X,Y))/\Hom_{\V_0}(I, \mathcal{P}_{\V}(X,Y))\\
& = \Hom_{\V_0}(I,    \underline{[\A, \V]}(X,Y)).
\end{align*}
 \end{proof}

 \subsection{Remark.}\label{remark_injective_stable} If $[\A,\V]$ exists,   that $\V_0$ is quasi-abelian, that the unit object $I$ is $\mathcal{E}_{\textrm{max}}$-projective in $\V_0$, and that $\V_0$ admits  $  \left| {\rm N}^+(A) \right |$-fold products for every $A\in \A$, then by analogous arguments, one can construct  a $\V$-subfunctor $\mathcal{I}_{\V}(-,-)$ of $[\A,\V](-,-)$, and a $\V$-enriched category $\overline{[\A, \V]}$ such that  $(\overline{[\A, \V]})_0=\St_{ \mathcal{I}}([\A, \V]_0)$.


\section*{Funding}
The author  was  partly supported by grant 22004/PI/22 funded by Fundaci\'on S\'eneca, Agencia de Ciencia y Tecnolog\'ia de la Regi\'on de Murcia and by grant PID2024-155576NB-I00 funded by MICIU/ AEI/10.13039/501100011033/FEDER, UE.

\bibliographystyle{plain}

\end{document}